\documentclass[12pt,a4paper]{article}
\usepackage{amsmath,amssymb,amsthm}
\usepackage{color,graphics}

\makeatletter

\newtheorem{theorem}{Theorem}[section]
\newtheorem{lemma}[theorem]{Lemma}
\newtheorem{corollary}[theorem]{Corollary}
\newtheorem{remark}[theorem]{Remark}

\theoremstyle{definition}
\newtheorem{de}[theorem]{Definition}
\newtheorem{example}[theorem]{Example}

\makeatother

\def\e{\begin{equation}}
\def\ee{\end{equation}}
\def\1{\frac{1}{2}}

\def\Id{\mathop{\mathrm{Id}}\nolimits}
\def\im{\mathop{\mathrm{Im}}\nolimits}
\def\Im{\mathop{\mathrm{Im}}\nolimits}
\def\Ker{\mathop{\mathrm{Ker}}\nolimits}
\def\tr{{\rm tr}\,}
\def\FF{{\mathbb F}}
\def\NN{{\mathbb N}}
\def\ZZ{{\mathbb Z}}
\def\CC{{\mathbb C}}
\def\IR{{\mathbb R}}

\def\rk{{\rm rk \,}}
\def\ch{{\rm char}\, }
\def\diag{{\rm diag}\,}
\def\Ss{\mathfrak{S}_{\mathfrak{p}}}

\def\Sp{\mathop{\mathrm{Sp}}}
\def\Lin{\mathop{\rm Lin}}
\def\DPhi{\mathfrak{D}}
\def\anyargument{\llcorner\!\!\!\lrcorner}

\begin{document}\showhyphens{further}
\title{Preserving zeros of  a
polynomial\thanks{The research of the first author was supported  by
the grant from Russian Ministry of Science and Technology,
MK-1417.2005.1. The research of the second author was supported by
grants from the Ministry of Science of Slovenia.}}
\author{A.\ E.\ Guterman$^{\rm a}$ \\
\parbox[b]{13cm}{\center\footnotesize Faculty of Algebra, Department of Mathematics and
Mechanics, Moscow State University, GSP-2, 119992 Moscow, Russia.} \and
B.\ Kuzma$^{\rm b}$\\
\parbox[b]{13cm}{\center\footnotesize ${}^1$Institute of Mathematics, Physics, and
Mechanics,
Jadranska~19, 1000~Ljubljana, Slovenia.\\
${}^2$University of Primorska, Faculty of Education, Cankarjeva 5, 6000
Koper, Slovenia.}}
\date{}
\maketitle \footnotetext[1]{E-mail: \texttt{guterman@list.ru}}
\footnotetext[2]{E-mail: \texttt{bojan.kuzma@pef.upr.si}}
\begin{abstract}
We study non-linear surjective mappings on subsets of~${\cal
M}_n(\FF)$, which preserve the zeros of some fixed polynomials in
noncommuting variables.


\end{abstract}

\section{Introduction}

The theory of transformations preserving different properties and
invariants dates back to the works by Frobenius, \cite{Fr}, Schur,
\cite{Sch}, and Dieudonn\'e, \cite{Die}, and is an intensively
developing part of algebra nowadays. A characterization of maps
preserving zeros of polynomials plays a central role in this area. The
most known results of this type were devoted to the characterization of
linear maps on matrix algebras, preserving the following polynomials:
$\mathfrak{p}(x)=x^k$, which correspond to nilpotent matrices,
see~\cite{BPW}, $\mathfrak{p}(x)=x^2-x$, whose zeros are idempotent
matrices, see~\cite{BP,CL,BS1}, potent matrices satisfying
$\mathfrak{p}(x)=x^k-x$, see~\cite{BS}, and matrices of finite order,
i.e., the zeros of $\mathfrak{p}(x)=x^k-1$, see~\cite{PD,GLS}. In 1980
Howard~\cite{howard} proved the general classification theorem for
bijective linear transformations on matrix algebras, preserving zeros
of a polynomial in one variable with at least two distinct zeros. This
result, together with the main theorem from~\cite{BPW}, provides the
complete characterization of bijective linear maps on matrix algebras
over algebraically closed fields, preserving zeros of polynomials in
one variable. Later in \cite{Li_Pierce} Li and Pierce investigated the
possibility to remove the invertibility assumption from Howard's
theorem and proved some related results.

In parallel since 1976, a question of characterizing linear
transformations preserving zeros of multilinear polynomials in several
noncommuting variables was considered. In particular, Watkins,
\cite{Watk}, characterized bijective linear transformations preserving
commutativity, i.e., zeros of the polynomial $\mathfrak{p}(x,y)=xy-yx$,
in \cite{Wong,Wong1} Wong classified operators preserving zero
products, i.e., zeros of $\mathfrak{p}(x,y)=xy$. During the last years
there was a big interest to this question, see
\cite{beasly_guterman_Lee_Song,beidar_lin,bresar,chan_li_sze,chebotar_fong_lee,fosner,petek_t,semrl:commutativity,zhao_hou}
and references therein. Also several related topics were intensively
investigated. Namely, additive transformations $\Phi$ on prime rings
satisfying the stronger condition $\mathfrak{p}(\Phi(x_1),\ldots,
\Phi(x_k))=\Phi(\mathfrak{p}(x_1,\ldots,x_k))$ for all $x_1,\ldots,
x_m$, were studied in~\cite{Beidar_Fong}. Linear maps that preserve
operators on infinitely dimensional algebras annihilated by a
polynomial in one variable, were treated
in~\cite{Hou_Hou,Semrl:Operator_Theory}.  Bijective linear operators on
the matrix algebra preserving zeros of the involutory polynomial
$\mathfrak{p}(x,y)=xy-yx^*$ are classified in \cite{chebotar_fong_lee}.
{\em Additive\/} surjections on certain classes of algebras, preserving
zeros of $\mathfrak{p}(x)=x^2$ or preserving zeros of the Jordan
polynomial $\mathfrak{p}(x,y)=xy+yx$, were characterized
in~\cite[Theorem~4.1]{Chebotar_Ke_Lee}, \cite[Lemma 2.3]{hou_zhao}, and
in~\cite{zhao_hou}.

In spite of the constant interest to the characterization problems for
linear transformations, preserving zeros of matrix polynomials, there
were no general answers in the multivariable case, like the Howard's
theorem concerning linear preservers of zeros of polynomials in one
variable. In particular, all results were proved for some specific
polynomials.

In \cite{chebotar_fong_lee} Chebotar, Fong and Lee posed the question
about the general form of linear preservers for zeros of multivariable
polynomials explicitly as an open problem:
\medskip

{\bf Problem}. \cite[Problem 1.1]{chebotar_fong_lee} Let
$\mathfrak{p}(x_1,\dots,x_k)$ be a polynomial over a field $\FF$ in
noncommuting indeterminates $x_1,\ldots, x_k$ of the degree $\deg
\mathfrak{p}>1$, and $\Phi: {\cal M}_n (\FF) \to {\cal M}_n(\FF)$ be a
linear map on the matrix algebra ${\cal M}_n (\FF)$. Suppose that
$\mathfrak{p}(A_1,\ldots, A_k )=0$ implies
$\mathfrak{p}(\Phi(A_1),\ldots , \Phi(A_k))=0$. Is it possible to
describe such $\Phi$?
\medskip

The authors of \cite{chebotar_fong_lee} have conjectured that if the
size of matrices, $n$, is big enough comparing with the number of
variables, $k$, and $\Phi$ is linear and bijective, then~$\Phi$ is a
sum of a scalar multiple of Jordan homomorphism and a transformation
that maps into the center of algebra.

The present paper is devoted to the solution of the above mentioned
problem for certain sufficiently large classes of polynomials of a
general type. Several remarks are in order. Firstly, our results are
non-linear in nature, we even do not assume additivity of a
transformation under consideration in advance. Secondly, the
transformation is not necessary required to be invertible and we only
assume that it is surjective. In addition we found some conditions
which replace the surjectivity assumption and also provide the examples
showing that the assumptions on the transformations, we have posed, are
indispensable. Moreover, the developed technique is characteristic free
and allows us to work without  restrictions on the number of variables
of a polynomial. This is done by the exclusion of polynomials which do
not provide sufficient restrictions on the transformation to be
classified. For example, this is the case with polynomial identities of
${\cal M}_n(\FF)$. Say, the polynomial
$\mathfrak{p}(x_1,\dots,x_{2n}):=\sum_{ \sigma\in{\cal
S}_{2n}}sign(\sigma)x_{\sigma(1)}\cdots x_{\sigma(2n)}$ is an identity
on ${\cal M}_n(\FF)$ by the famous Amitsur-Levitzki's theorem,
see~\cite{Amitsur_Levitzki}. This polynomial  vanishes on the whole
matrix algebra and therefore it gives no condition on~$\Phi$. Therefore
we divide our considerations in two parts. Firstly we consider the {\em
generic case\/}, where the sum of coefficients of a multilinear
polynomial ${\mathfrak p}$ is non-zero, and thus ${\mathfrak p}$ can
not be an identity in ${\cal M}_n(\FF)$. Then we investigate the {\em
derogatory case\/}, where the sum of coefficients of ${\mathfrak p}$ is
zero and polynomial identities may appear.

\medskip

Throughout, $n\ge 3$ will be an integer and ${\cal M}_n({\FF})$ will be
the algebra of $n\times n$-matrices over an arbitrary field ${\FF}$.
Let $E_{ij}$ be its standard basis. Let $\mathrm{GL}_n({\FF}) \subset
{\cal M}_n({\FF})$ denote the group of invertible matrices, with
identity~$\Id$. Let~${\cal I}^1\subseteq {\cal M}_n(\FF)$ be the set of
all rank-one idempotents. Given a field
homomorphism~$\varphi:\FF\to\FF$ ({\small i.e., an additive and
multiplicative function on $\FF$}), we let~$X^\varphi$ be a matrix,
obtained from~$X$ by applying~$\varphi$ entry-wise. In addition,
let~$X^{\tr}$ be the transposed matrix of~$X$. Matrices $P,Q\in {\cal
M}_n(\FF)$ are called {\em orthogonal\/} to each other if $PQ=QP=0$.

Lastly, let~${\cal S}_k$ be the set of all permutations of the set
$\{1,\dots,k\}$.\bigskip

\medskip

\begin{de} Let $k\ge 2$.  We say that a matrix $k$-tuple $(A_1,\ldots, A_k)$ {\em
is a zero of a homogeneous multilinear polynomial}
$${\mathfrak{p}}(x_1,\dots,x_k):=\sum_{\sigma\in{\cal S}_k}
\alpha_\sigma x_{\sigma(1)}\cdots x_{\sigma(k)} $$ if
$$\mathfrak{p}(A_1,\dots, A_k ):=\sum_{\sigma\in{\cal S}_k} \alpha_\sigma
A_{\sigma(1)} \cdots A_{\sigma(k)}=0.$$
The set of all such $k$-tuples will be denoted by~$\Ss\subseteq {\cal
M}_n(\FF)\times\dots\times {\cal M}_n(\FF)$.
\end{de}
\begin{de}
Suppose~$\mathfrak{p}_1,\mathfrak{p}_2$ are  two homogenous multilinear
polynomials. A transformation $\Phi:{\cal M}_n(\FF)\to {\cal M}_n(\FF)$
{\em maps the zeros of~$\mathfrak{p}_1$ to the zeros of
$\mathfrak{p}_2$} whenever the implication $\mathfrak{p}_1(A_1,\dots,
A_k )=0\Longrightarrow\mathfrak{p}_2(\Phi(A_1),\dots, \Phi(A_k) )=0$
holds ({\small equivalently, whenever $\Phi({\Ss}_1)\subseteq
{\Ss}_2$}). If~$\mathfrak{p}_1=\mathfrak{p}_2=:\mathfrak{p}$ then
$\Phi$ {\em preserves the zeros of~$\mathfrak{p}$}. In addition, if
$\mathfrak{p}(x,y)=xy-yx$ then $\Phi$  {\em preserves commutativity}.
\end{de}
\begin{de}A transformation $\Phi:{\cal M}_n(\FF)\to {\cal M}_n(\FF)$
{\em strongly maps the zeros of~$\mathfrak{p}_1$ to the zeros
of~$\mathfrak{p}_2$} whenever the equivalence
$\mathfrak{p}_1(A_1,\dots, A_k
)=0\Longleftrightarrow\mathfrak{p}_2(\Phi(A_1),\dots, \Phi(A_k) )=0$
holds. If $\mathfrak{p}_1=\mathfrak{p}_2=:\mathfrak{p}$ then $\Phi$
{\em strongly preserves the zeros of~$\mathfrak{p}$}.
\end{de}

\medskip

Our paper is organized as follows.

In Section~\ref{S4} we study the mappings between certain matrix
subspaces, including the map from the whole matrix algebra to itself,
which strongly preserve zeros of homogeneous multilinear polynomials
with nonzero sum of coefficients.

In Section~\ref{S3} we study the transformations that map zeros of a
homogeneous multilinear polynomials of arbitrary many variables with
zero sum of coefficients to zeros of another such polynomial. To avoid
the obstructions which come from the polynomial identities of matrix
algebra, it is necessary to restrict the set of permutations. The
general problem is then reduced to the already well-studied
commutativity preservers,
see~\cite{semrl:commutativity,fosner,bresar,petek_t,beidar_lin} for
their characterization. We remark that some ideas that we are using in
this section came  from our recent paper~\cite{alieve_guterman_kuzma}.

Section~\ref{S5} contains a number of examples showing  that our
assumption are indispensable without posing some additional conditions
on $\Phi$ or {on $\mathfrak{p}(x_1,\ldots,x_k)$}.

\section{\label{S4} Polynomials with non-vanishing sums of\\ coefficients}

In the present section we investigate surjective mappings on certain
matrix subspaces, in particular on the whole matrix space, that {\em
strongly} preserve zeros of a polynomial with non-vanishing sum of
coefficients.

We also refer to Chan, Li, and Sze~\cite{chan_li_sze}, where the zeros
of a polynomial $\mathfrak{p}(x,y):=xy$ were considered, and, similarly
to our results below, the nice structure was obtained solely on
matrices of rank-one. We will show in the last section that in a way
our results cannot be further improved, without imposing additional
hypothesis, say additivity of~$\Phi$.

However, if~$\Phi$ strongly preserves the zeros of a polynomial with at
least three variables, then we were able to deduce a structural result
holding {\em for all matrices} from the defining set of~$\Phi$.\bigskip

We now list the main results of the present section.
Let~$\DPhi_1\subseteq {\cal M}_n(\FF)$ and $\DPhi_2\subseteq {\cal
M}_n(\FF)$ be subsets {\em that contain all matrices of rank-one and
all idempotents of rank $n-1$\/}. Also in this section we assume that a
homogeneous multilinear polynomial
\begin{equation}\label{eq:poly-with-nonvanishing-sum}%
\mathfrak{p}(x_1,\dots ,x_k):=\sum_{\sigma\in{\cal S}_k} \alpha_\sigma
x_{\sigma(1)}\cdots x_{\sigma(k)};\qquad (\alpha_\sigma\in\FF)
\end{equation}
satisfies $\sum \alpha_\sigma\not=0$.\medskip

 The most general form of the main result of this section is the following:
\begin{theorem}\label{thm:main-nonvanishing-sum}%
Let~$\FF$ be an arbitrary algebraically closed field. Assume~$n\ge 3$
and $k\ge 2$. If a surjection $\Phi:\DPhi_1\to\DPhi_2$ strongly
preserves the zeros of~$\mathfrak{p}(x_1, \ldots,x_k)$ then there
exists a field isomorphism $\varphi:\FF\to\FF$, a function
$\gamma:\DPhi_1\backslash\{0\}\to\FF^\ast:=\FF\backslash\{0\}$, and an
invertible matrix~$T$ such that
\begin{itemize}
\item[{\rm (i)}]
$\Phi(A) = \gamma(A) \,
 {T}A^{\varphi}{T}^{-1}$ for all rank-one
matrices~$A$, or
\item[{\rm (ii)}]
$\Phi(A) = \gamma(A) \,
 {T}\left(A^{\varphi}\right) ^{\tr}{T}^{-1}$
for all rank-one matrices~$A$.
\end{itemize}
\end{theorem}
\begin{remark} {The {converse to Theorem~\ref{thm:main-nonvanishing-sum}}
does not hold without imposing additional constraints on
isomorphism~$\varphi$. Namely}, there are transformations  of types (i)
and (ii) which do not preserve the zeros of~$\mathfrak{p}$. We refer
the reader to the last section for examples.
\end{remark}

When  $\ch\FF\not=2$, and the dimension of matrices is $n\ge 4$, and
the polynomial has at least three variables we have a nice structural
result holding {\em for all matrices} from the defining set~$\DPhi_1$
of~$\Phi$. In particular, for all matrices if $\DPhi_1={\cal
M}_n(\FF)$. We merely add scalar matrices to conclusions (i) and (ii)
of Theorem~\ref{thm:main-nonvanishing-sum}:

\begin{corollary} \label{cor:TAllMatr}%
Under the assumptions of  Theorem~\ref{thm:main-nonvanishing-sum},
assume further $\ch\FF\not=2$, and $n\ge 4$, and  \mbox{$k\ge 3$}.
Then,
\begin{itemize}
\item[{\rm (i)}]
$\Phi(A) = \gamma(A)\, {T}A^{\varphi}{T}^{-1}+\mu(A)\,\Id$ for all
$A\in {\DPhi_1}$, or
\item[{\rm (ii)}]
$\Phi(A) = \gamma(A) \, {T}\left(A^{\varphi}\right)
^{\tr}{T}^{-1}+\mu(A)\,\Id$ for all $A\in {\DPhi_1}$.
\end{itemize}
\end{corollary}

The situation is completely different when $k=2$,  see
Example~\ref{exa:counter-axa-for-nonvanishing} below. However,
if~$\mathfrak{p}(x,y):=xy+yx$ is a  polynomial of Jordan multiplication
we can still get the characterization for some special
matrices~$A\in\DPhi_1$.

\begin{corollary}\label{cor:corolary_to_Jordan_product}%
Under the assumptions of  Theorem~\ref{thm:main-nonvanishing-sum},
assume further $\mathfrak{p}(x,y):=xy+yx$ and~$\ch\FF\not=2$. Then, the
conclusions (i) and (ii) also hold for diagonalizable~$A\in\DPhi_1$,
with the spectrum, $\Sp(A)=\{\lambda, -\lambda\}$.
\end{corollary}

Moreover, we may remove the surjectively assumption from the
Theorem~\ref{thm:main-nonvanishing-sum}, at least for some polynomials:
\begin{corollary}\label{cor:nonsurjective}%
Under the assumptions of  Theorem~\ref{thm:main-nonvanishing-sum},
suppose that a possibly nonsurjective~$\Phi:\DPhi_1\to\DPhi_2$ strongly
preserves the zeros of polynomial~$\mathfrak{p}$, defined in
Eq.~(\ref{eq:poly-with-nonvanishing-sum}), but such that the matrix
$$\mathrm{Cof}\,(\mathfrak{p}):=\begin{pmatrix}
\sum_{\sigma(1)=1}\alpha_\sigma &\sum_{\sigma(1)=2} \alpha_\sigma&\dots &\sum_{\sigma(1)=k}\alpha_\sigma\\
\sum_{\sigma(2)=1}\alpha_\sigma &\sum_{\sigma(2)=2} \alpha_\sigma&\dots &\sum_{\sigma(2)=k}\alpha_\sigma\\
\vdots & \vdots &\ddots &\vdots\\
\sum_{\sigma(k)=1}\alpha_\sigma &\sum_{\sigma(k)=2} \alpha_\sigma&\dots
&\sum_{\sigma(k)=k}\alpha_\sigma
\end{pmatrix}\in{\cal M}_k(\FF)$$
is invertible. Then, the conclusions (i) and (ii) remain valid, with
the exception that a field homomorphism $\varphi:\FF\to\FF$ might be
nonsurjective.
\end{corollary}
\subsection{The proof of Theorem~\ref{thm:main-nonvanishing-sum}}

We first rewrite~$\mathfrak{p}$. Namely, its coefficients satisfy
$$\sum\alpha_\sigma=\sum_{\sigma(1)=1}\alpha_\sigma+\sum_{\sigma(2)=1}
\alpha_\sigma+\dots+ \sum_{\sigma(k)=1}\alpha_\sigma,$$
 and since the left-hand side is nonzero, so must be at least one of
the summands on the right. Say,
$0\not=\sum_{\sigma(\imath_0)=1}\alpha_\sigma$. By
dividing~$\mathfrak{p}$ we may assume
$\sum_{\sigma(\imath_0)=1}\alpha_\sigma=1$. Moreover, we may also
assume $\imath_0=1$. Otherwise we would regard the polynomial
$\mathfrak{p}'(x_1,\dots,x_k):=\mathfrak{p}(x_{\tau(1)},\dots,x_{\tau(k)})$
for  permutation~$\tau:=(1,\imath_0)$. Obviously,~$\Phi$ would still
strongly preserve the zeros of~$\mathfrak{p}'$. We can now
rewrite~$\mathfrak{p}$  into the form
\begin{equation}\label{eq:p}%
\mathfrak{p}(x_1,\dots,x_k):=\sum_{\sigma(1)=1}\alpha_\sigma
x_{\sigma(1)}\dots x_{\sigma(k)} +\sum_{\sigma(1)\not=1}\alpha_\sigma
x_{\sigma(1)}\dots x_{\sigma(k)};
\end{equation}
where $\sum_{\sigma(1)=1}\alpha_\sigma=1$, and where
$\xi:=\sum_{\sigma(1)\not=1}\alpha_\sigma\not=-1$.  Similarly, there
exists $\jmath=\jmath_0$ such that
$$\xi_{\jmath_0k}:=\sum_{\sigma(\jmath_0)=k} \alpha_\sigma\not=0,$$ for
otherwise $0=\sum_{\sigma(k)=k} \alpha_\sigma+\sum_{\sigma(k-1)=k}
\alpha_\sigma+\dots +\sum_{\sigma(1)=k}
\alpha_\sigma=\sum\alpha_\sigma,$ a contradiction! In the sequel, we
will always use the equivalent form~(\ref{eq:p}) of
polynomial~$\mathfrak{p}$.
\bigskip

{It will be beneficial for our further considerations to associate with
each matrix~$A$ the  two sets: $\Omega_{\bullet A}$ and
$\Omega_{A\bullet A}$, defined via the polynomial~$\mathfrak{p}$ by}
\begin{align}\label{eq:def_Omega_A}%
\Omega_{\bullet A}&:=\{X\in{\cal
M}_n(\FF);\;\;\mathfrak{p}(X,A,\dots,A)=0\},\\
\Omega_{A\bullet A }&:=\{X\in{\cal
M}_n(\FF);\;\;\mathfrak{p}(A,\dots,A,\fbox{$X$}_{\jmath_0},A,\dots,A)=0\},
\label{eq:def_Omega_A_1}
\end{align}
in the last equation,~$X$ is at the ${\jmath_0}$ position. Note that
$\Omega_{A\bullet A }=\Omega_{\bullet A}$ if~$\jmath_0=1$ ({\small  say
for polynomial
$\mathfrak{p}(x_1,x_2):=x_1x_2-x_2x_1$}).\typeout{MyWarning: CHANGE on
22.2.07 input line \the\inputlineno.} Clearly,~$\Omega_{\bullet A}$ and
$\Omega_{A\bullet A }$ are both vector subspaces of~${\cal M}_n(\FF)$.
Moreover,
 we can rewrite the condition
$X\in\Omega_{\bullet A}$, respectively, $X\in\Omega_{A\bullet A}$, as
\begin{align}\label{eq:pom3}%
    XA^{k-1}+(\beta_2\,AXA^{k-2}+\dots+\beta_{k}\,A^{k-1}X)&=0,
    \quad\Bigl(\beta_i:=\sum_{\sigma(1)=i}\alpha_\sigma\Bigr)\\
 \intertext{respectively,}
    \xi_{\jmath_0 k}\,A^{k-1}X+(\tilde{\beta}_{k-1}\,A^{k-2}XA+\dots+
        \tilde{\beta}_{1}\,XA^{k-1})&=0,\quad\Bigl(\tilde{\beta}_i:=
    \sum_{\sigma(\jmath_0)=i}\alpha_\sigma\Bigr).\label{eq:pom3.1}%
\end{align}
In particular,~$\Omega_{\bullet A}$ equals the null space of the
elementary operator {$\mathbf{T}_{\bullet A}:{\cal M}_n(\FF)\to{\cal
M}_n(\FF)$}, defined by $X\mapsto
XA^{k-1}+(\beta_2\,AXA^{k-2}+\dots+\beta_{k}\,A^{k-1}X)$.

Similarly, $\Omega_{A\bullet A}$ equals the null space of of the
elementary operator {$\mathbf{T}_{A\bullet A} $}, defined by $X\mapsto
A^{k-1}X+(\hat{\beta}_{k-1}\,A^{k-2}XA+\dots+\hat{\beta}_{1}\,XA^{k-1})$,
where $\hat{\beta}_i:=\tilde{\beta}_i/\xi_{\jmath_0k}$.\bigskip

We proceed with a series of lemmas. The first lemma allows us to
compute the spectrum of elementary operators, in particular, of the
operators $\mathbf{T}_{\bullet A}$ and $\mathbf{T}_{A \bullet A}$. We
present the easy proof for the sake of convenience.
\begin{lemma}\label{lem:Sylvester--Rosenblum}%
Let $\FF$ be an arbitrary field. Let $L\in{\cal M}_m(\FF)$ and $M\in
{\cal M}_n(\FF)$ be matrices with $\Sp(L)=\{\lambda\}$ and
$\Sp(M)=\{\mu\}$. Then, the spectrum of elementary operator
$\mathbf{T}_{L,M}:{\cal M}_{m\times n}(\FF)\to{\cal M}_{m\times
n}(\FF)$, defined by $ X\mapsto \beta_t L^tX+\beta_{t-1}
L^{t-1}XM+\dots+\beta_0 XM^t$ is a singleton:
$\Sp(\mathbf{T}_{L,M})=\{\beta_t \lambda^t+\beta_{t-1}
\lambda^{t-1}\mu+\dots+\beta_0\mu^t\}$.
\end{lemma}
\begin{proof}
Let us consider the decomposition~$L=S_L\hat{L}S_L^{-1}$,
where~$\hat{L}$ is the upper-triangular Jordan form of~$L$, and $M=S_M
\breve{M}S_M^{-1}$, where $ \breve{M}$ is lower-triangular Jordan form
of~$M$. It is well-known that the matrix representation of $X\mapsto
L^{t-k}XM^k$, relative to basis
$$E_{11},\dots,E_{n1},E_{12},\dots,E_{n2},\dots,E_{1n},\dots,E_{nn}$$
equals the  tensor (Kronecker) product $(M^k)^{\tr}\otimes L^{t-k}$.
Moreover, the matrices $(M^k)^{\tr}\otimes L^{t-k}$; $(k=0,1,\dots)$
are simultaneously similar to~$(\breve{M}^k)^{\tr}\otimes
\hat{L}^{t-k}$ via similarity~$(S_M^{-1})^{\tr}\otimes S_L$.

Hence, the spectrum of $(M^k)^{\tr}\otimes L^{t-k}$ equals the spectrum
of $(\breve{M}^k)^{\tr}\otimes \hat{L}^{t-k}$. Now,
$(\breve{M}^k)^{\tr}$, as well as $\hat{L}^{t-k}$ are both
upper-triangular matrices, with $\mu^k$, respectively, $\lambda^{t-k}$
 on the diagonal. Hence, their tensor product remains
upper-triangular, with~$\lambda^{t-k}\mu^k$  on the diagonal.
Consequently, $\mathbf{T}_{L,M}=\sum \beta_k (M^k)^{\tr}\otimes
L^{t-k}$ is similar to the upper-triangular matrix $\sum \beta_k
(\breve{M}^k)^{\tr}\otimes \hat{L}^{t-k}$, with the number $\sum\beta_k
\lambda^{t-k} \mu^k$ on main diagonal. This number is, hence, the only
eigenvalue of $\mathbf{T}_{L,M}$.
\end{proof}
We remark that over the field of complex numbers this fact follows from
the results of  Lumer and Rosenblum~\cite{lumer_rosenblum} and
Curto~\cite{curto}, where it was proven in a different way  for linear
operators on Hilbert spaces, possibly infinite dimensional. We also
remark that in the case $\mathbf{T}_{\bullet A}=A^{k-1}X+XA^{k-1}$ a
short proof of the corresponding result for Hilbert spaces is presented
by Bhatia and Rosenthal in~\cite[p.~2]{bhatia_rosenthal} and some
further properties of the spectrum can be found in~\cite{curto}, some
additional properties of $\mathbf{T}_{\bullet A}$ are investigated by
Chuai and Tian in~\cite{chuai_tian}.

\begin{lemma}\label{lem:A-circ-B=0<=>AB=0=BA}%
Let $\mu,\nu,\mu',\nu'\in\FF$ be scalars, and let
 $1+\mu+\nu\not=0$.
If  $$PX+\mu\,PXP+\nu XP=0=XP+\mu'\,PXP+\nu' PX$$
 holds for some idempotent~$P$ then $PX=0=XP$.
\end{lemma}
\begin{proof}
Postmultiply the equation $PX+\mu\,PXP+\nu XP=0$ with $P$ and subtract.
We derive $PX=PXP$. Likewise $PXP=XP$, from the second equation. Using
$PX=PXP=XP$ again in the first equation, we get $(1+\mu+\nu)PXP=0$, so
$XP=PXP=PX=0$.
\end{proof}
\begin{lemma}\label{lem:Phi(X)=0<==>X=0}%
 Suppose $\DPhi\subseteq{  M}_n(\FF)$ contains all matrices of rank-one.
 Then, $A\in {\cal  M}_n(\FF)$ is a zero matrix if and only if
$\mathfrak{p}(A,X,\dots,X)=0$ for each~$X\in\DPhi$.
\end{lemma}
\begin{proof}
We prove only the nontrivial implication. Indeed,
substituting~$X:=E_{ii}$ we have
\begin{align*}
0&=\mathfrak{p}(A,E_{ii},\dots, E_{ii})=\sum_{\sigma(1)=1}\alpha_\sigma
AE_{ii}+ \sum_{\sigma(1)\not\in\{1,k\}}\alpha_\sigma
E_{ii}AE_{ii}+\sum_{\sigma(1)=k}\alpha_\sigma E_{ii}A\\
&= AE_{ii}+\alpha E_{ii}AE_{ii}+\beta_kE_{ii}A
\end{align*}
 for $i=1,\ldots, n$.
Premultiply with idempotent~$E_{ii}$, and compare the two equations. We
get $A E_{ii}=E_{ii}AE_{ii}$. Hence, $0=AE_{ii}+(\alpha
AE_{ii}+\beta_{k} E_{ii}A)$. We may sum-up these $n$ equations to get
$0=A\Id+(\alpha A\Id+\beta_{k}\Id A)=(1+\alpha +\beta_{k})A$, so~$A=0$.
\end{proof}
\begin{corollary}\label{cor:Phi-preserves-0}%
Let $\Phi$ satisfy conditions of
Theorem~\ref{thm:main-nonvanishing-sum}.  Then, $0\in\DPhi_1$ if and
only if $0\in\DPhi_2$. Moreover, $\Phi(A)=0$ if and only if~$A=0$.
\end{corollary}
\begin{proof}
If~$A\not=0$ then, by Lemma~\ref{lem:Phi(X)=0<==>X=0}, there exists
some $X\in\DPhi_1$ with $\mathfrak{p}(A,X,\dots,X)\not=0$.
Consequently, $\mathfrak{p}(\Phi(A),\Phi(X),\dots,\Phi(X))\not=0$, and
so $\Phi(A)\not=0$. Hence, if~$0\not\in\DPhi_1$ then
$0\notin\im\Phi=\DPhi_2$. By surjectivity we likewise see
that~$\Phi(A)\not=0$ implies~$A\not=0$.
\end{proof}

We now characterize when two rank-one nilpotents are scalar multiple of
each other in terms of~$\Omega_{\bullet A}\cap\Omega_{A\bullet A}$,
i.e., in terms of the zeros of polynomial~$\mathfrak{p}$.

\begin{lemma}\label{lem:charact_(nilpot-of-rk1)_with_idempot}%
Let $n\ge 3$, let $\varphi:\FF\to\FF$ be a nonzero field homomorphism,
and let $N_1,N_2\in{\cal M}_n(\FF)$ be rank-one nilpotents. Assume that
the following condition (i) is  satisfied:
\begin{itemize}
\item[(i)] $N_1\in\Omega_{\bullet P}\cap\Omega_{P\bullet P}
     \Longleftrightarrow N_2\in\Omega_{\bullet P^\varphi}\cap
    \Omega_{P^\varphi\bullet P^\varphi}$ holds for     every rank-one
    idempotent~$P$.
\end{itemize}
Then $N_2=\lambda N_1^\varphi$ for some nonzero scalar~$\lambda\in
\FF$.
\end{lemma}
\begin{proof}
Pick a similarity~$S$ such that $N_1 =SE_{12}S^{-1}$. Then, the
rank-one idempotents $P_3:=SE_{33}S^{-1},\dots, P_n:=S E_{nn}S^{-1}$,
and $P_{n+1}:=S(E_{n2}+E_{nn})S^{-1}$ are orthogonal to~$N_1$. Hence,
$N_1\in\Omega_{\bullet P_i}\cap\Omega_{P_i\bullet P_i}$, which by (i)
implies $N_2\in\Omega_{\bullet P_i^\varphi}\cap
\Omega_{P_i^\varphi\bullet P_i^\varphi}$ for $i=3,\ldots, (n+1)$. Using
the equivalent expressions~(\ref{eq:pom3})--(\ref{eq:pom3.1}),  we can
easily rewrite this into
$$\begin{aligned}
0&=\bigl(  (S^{-1})^\varphi N_2S^\varphi\bigr)\circ_\bullet E_{33} &
0&=\bigl( (S^{-1})^\varphi
   N_2S^\varphi\bigr)\mathop{_\bullet\circ}E_{33}\\
\multispan4{\dotfill}\\
 0&=\bigl((S^{-1})^\varphi N_2S^\varphi)\circ_\bullet E_{nn}
  &
0&=\bigl((S^{-1})^\varphi N_2S^\varphi)\mathop{_\bullet\circ}E_{nn}\\
0&=\bigl( (  S^{-1})^\varphi N_2S^\varphi\bigr)\circ_\bullet
(E_{n2}+E_{nn}) & 0&=\bigl( (  S^{-1})^\varphi
N_2S^\varphi\bigr)\mathop{_\bullet\circ}(E_{n2}+E_{nn}),
\end{aligned}$$
where $X\circ_\bullet P:=XP+\beta PXP+\beta_k PX$, and where
$X\mathop{_\bullet\circ}P:=PX+\hat{\beta}\,PXP+ \hat{\beta}_{1}\,XP$
for scalars $\beta:=(\beta_{2}+\dots+\beta_{k-1})$,
$\hat{\beta}:=(\tilde{\beta}_{k-1}+\dots+\tilde{\beta}_{2})/\xi_{\jmath_0k}$,
and $\hat{\beta}_1:=\tilde{\beta}_1/\xi_{\jmath_0k}$.

By Lemma~\ref{lem:A-circ-B=0<=>AB=0=BA}, the above equations imply
orthogonality between the nilpotent $(S^{-1})^\varphi N_2S^\varphi$ and
idempotents $E_{33},\dots,E_{nn}, (E_{n2}+E_{nn})$. More precisely, the
first $(n-2)$ equalities give  that $(S^{-1})^\varphi N_2S^\varphi$ can
be nonzero only in the upper-left $2\times 2$ block, while the last one
further yields $(S^{-1})^\varphi N_2S^\varphi=\varsigma E_{11}+\lambda
E_{12}$. Since~$N_2$ is nilpotent, $\varsigma=0$. Thus $N_2=\lambda
S^\varphi E_{12}(S^{-1})^\varphi= \lambda S^\varphi
E_{12}^\varphi(S^{-1})^\varphi =\lambda N_1^\varphi$, as desired.
\end{proof}
Similarly we can prove the following:
\begin{lemma}\label{lem:charact_(nilpot-of-rk1)_with_idempot'}%
Under the assumptions of
Lemma~\ref{lem:charact_(nilpot-of-rk1)_with_idempot}, suppose the
following condition (i') is satisfied:
\begin{itemize}
\item[(i')]$N_1\in\Omega_{\bullet P}\cap\Omega_{P\bullet P}
     \Longleftrightarrow N_2\in\Omega_{\bullet (P^\varphi)^{\tr}}\cap
    \Omega_{(P^\varphi)^{\tr}\bullet (P^\varphi)^{\tr}}$ holds for every
    rank-one idempotent~$P$.
\end{itemize}
Then $N_2=\lambda (N_1^\varphi)^{\tr}$ for some nonzero
scalar~$\lambda\in \FF$.
\end{lemma}
\begin{proof}
Similar to Lemma~\ref{lem:charact_(nilpot-of-rk1)_with_idempot}.
\end{proof}\bigskip

We next characterize scalar multiples of rank-one idempotents in terms
of~$\Omega_{\bullet A}\cap\Omega_{A\bullet A}$, i.e., in terms of  the
zeros of polynomial~$\mathfrak{p}$. This is a chief Lemma in the proof
of Theorem~\ref{thm:main-nonvanishing-sum}.
\begin{lemma}\label{lem:charact_of_minimal_idempotents}%
Fix~$A\in{\cal M}_n(\FF)$. Under the assumptions of
Theorem~\ref{thm:main-nonvanishing-sum}, {\em precisely one} of the
following three possibilities occurs for the set~$\Omega_{\bullet
A}\cap \Omega_{A\bullet A}$:
\begin{enumerate}
\renewcommand{\theenumi}{(\roman{enumi})}%
\item $\Omega_{\bullet A}\cap
\Omega_{A\bullet A}=\{0\}$.
\item $\Omega_{\bullet A}\cap
\Omega_{A\bullet A}$ contains a square-zero matrix of rank-one.
\item $\Omega_{\bullet A}\cap
\Omega_{A\bullet A}=\FF\,P$ where $P$ is a rank-one idempotent.
\end{enumerate}
\end{lemma}

\begin{proof}Obviously, the listed three possibilities are exclusive.
It, hence, remains to see that at least one of them does occur.
 We will rely on the fact that $X\in\Omega_{\bullet A}\cap \Omega_{A\bullet A}$ is
equivalent to Eq.~(\ref{eq:pom3}) and Eq.~(\ref{eq:pom3.1}),
simultaneously.

Now, with the help of similarity we may assume
$A=C_{n_1}(\lambda_1)\oplus\dots\oplus C_{n_r}(\lambda_r)$ is already
in its Jordan block-diagonal form, where~$\lambda_1,\dots,\lambda_r$
are pairwise distinct eigenvalues of~$A$, and an $n_i\times n_i$ matrix
$C_{n_i}(\lambda_i)$ is a sum of all Jordan blocks that correspond to
eigenvalue~$\lambda_i$. We may decompose~$X=\bigl(X_{ij} \bigr)_{1\le
i,j,\le r}$ accordingly. Then,
$A^{s}=C_{n_1}(\lambda_1)^{s}\oplus\dots\oplus C_{n_r}(\lambda_r)^{s}$,
for $s=1,2,\dots,n$, are also block-diagonal, so the $(i,j)$-th block
of Eqs.~(\ref{eq:pom3})--(\ref{eq:pom3.1}) read
\begin{align}\label{eq:pom4}%
X_{ij}C_{n_j}(\lambda_j)^{k-1}+
\bigl(\beta_2\,C_{n_i}(\lambda_i)X_{ij}C_{n_j}(\lambda_j)^{k-2}+\dots+
 \beta_{k}\,C_{n_i}(\lambda_i)^{k-1}X_{ij} \bigr)&=0
 \intertext{respectively,}
\xi_{\jmath_0k}C_{n_i}(\lambda_i)^{k-1}X_{ij}\!+\!
\bigl(\tilde{\beta}_{k-1}\,C_{n_i}(\lambda_i)^{k-2}X_{ij}C_{n_j}(\lambda_j)+
\!\cdots\!+
 \tilde{\beta}_{1}\,X_{ij}C_{n_j}(\lambda_j)^{k-1} \bigr)&\!=\!0\label{eq:pom4.1}%
 \end{align}
Obviously, $\Sp\bigl(C_{n_i}(\lambda_i)^{s}  \bigr)=\{\lambda_i^s\}$
for any integer~$s$. In view of Lemma~\ref{lem:Sylvester--Rosenblum} we
consequently introduce two polynomials
$$p_{\bullet
A}(\lambda,\mu):=\mu^{k-1}+(\beta_2\,\lambda\mu^{k-2}+\dots+\beta_{k}\,
\lambda^{k-1})\in\FF[\lambda,\mu]$$
as well as its counterpart
$$p_{A\bullet A}(\lambda,\mu):=\xi_{\jmath_0k}\lambda^{k-1}+
(\tilde{\beta}_{k-1}\,\lambda^{k-2}\mu+\dots+\tilde{\beta}_{1}\,\mu^{k-1}),$$
and  proceed in five steps:\medskip

\noindent {\bf Step 1}. Assume first that for no pair
$(\lambda_i,\lambda_j)\in\Sp(A)\times\Sp(A)$ we have simultaneously
$p_{\bullet
 A}(\lambda_i,\lambda_j)=0=p_{A\bullet A}(\lambda_i,\lambda_j)$.

In this case, we note that the left-hand sides of each of the
Eqs.~(\ref{eq:pom4})--(\ref{eq:pom4.1}) define an elementary operator.
Therefore, Lemma~\ref{lem:Sylvester--Rosenblum}, with
$L:=C_{n_i}(\lambda_i)$ and $M:=C_{n_j}(\lambda_j)$ implies that their
spectra are equal to $\{p_{\bullet A}(\lambda_i,\lambda_j)\}$, and
$\{p_{A\bullet A}(\lambda_i,\lambda_j)\}$, respectively. At least one
does not contain zero, and therefore, the  corresponding elementary
operator is invertible. The corresponding equation, on the other hand,
 implies that $X_{ij}=0$. Hence, all blocks of~$X$ are zero. This
clearly demonstrates $\Omega_{\bullet A}\cap \Omega_{A\bullet A
}=\{0\}$, and we have condition~(i) satisfied.\medskip

 \noindent{\bf Step 2}. Suppose next
$p_{\bullet A}(\lambda_i,\lambda_j)=0=p_{A\bullet
A}(\lambda_i,\lambda_j)$ for {\em distinct} $\lambda_i,\lambda_j$. For
simplicity, we assume $(i,j)=(1,2)$.

Here we consider the matrix~$X$ with all, but the $(1,2)$-th, blocks
zero. Then, by Eq.~(\ref{eq:pom4}), all blocks of
$\mathfrak{p}(X,A,\dots,A)$, but the $(1,2)$-th, are also zero. On the
other hand, its $(1,2)$ block equals to
\begin{align*}
\bigl(\mathfrak{p}(X,A,\dots,A)
\bigr)_{12}&=X_{12}C_{n_2}(\lambda_2)^{k-1}+\\
&\mbox{}+
\bigl(\beta_2\,C_{n_1}(\lambda_1)X_{12}C_{n_2}(\lambda_2)^{k-2}+\dots+
 \beta_{k}\,C_{n_1}(\lambda_1)^{k-1}X_{12} \bigr)
\end{align*}
 Now, write
$C_{n_1}(\lambda_1)=\lambda_1 \Id_{n_1}+N_1$ and
$C_{n_2}(\lambda_2)=\lambda_2 \Id_{n_2}+N_2$, for some upper-triangular
nilpotents~$N_1,N_2$. Next, consider an $n_1 \times n_2$ matrix
$\hat{X}_{12}$, with all entries, but the upper-right one, equal to
zero. It is easy to see that $\hat{X}_{12}N_2=0=N_1\hat{X}_{12}$, which
in turn, implies that the right side of the above equation simplifies
into $
\hat{X}_{12}\lambda_2^{k-1}+(\beta_2\,\lambda_1X_{12}\lambda_2^{k-2}+\dots+
 \beta_{k}\,\lambda_1^{k-1}X_{12}
 )=(\lambda_2^{k-1}+\beta_2\,\lambda_1\lambda_2^{k-2}+
 \dots+\beta_{k}\lambda_1^{k-1})\hat{X}_{12}=p_{\bullet
 A}(\lambda_1,\lambda_2)\hat{X}_{12}=0$.
Consequently, if a block-matrix $X\in {\cal M}_n(\FF)$ has its
$(1,2)$-th block equal to $\hat{X}_{12}$ while the other blocks are
zero then it is  a rank-one nilpotent, and
$\mathfrak{p}(X,A,\dots,A)=0$. {Similarly, by Eq.~(\ref{eq:pom4.1}), we
also infer
$\mathfrak{p}(A,\dots,A,\fbox{$X$}_{\jmath_0},A,\dots,A)=0$}.
Therefore, such~$X\in\Omega_{\bullet A}\cap \Omega_{ A\bullet A}$,
which guaranties the condition~(ii).\medskip

\noindent{\bf Step 3}.  Suppose  we are not under the conditions of
Step 2. We, thus, consider the case $p_{\bullet
 A}(\lambda_i,\lambda_i)=0=p_{A\bullet A}(\lambda_i,\lambda_i)$ for some $i$.

Clearly, $p_{\bullet
 A}(\lambda_i,\lambda_i)=  \lambda_i^{k-1}(1+\xi)$ where
$\xi:=\beta_2+\dots+\beta_{k}$. However, $\xi\not=-1$, so  $p_{\bullet
A}(\lambda_i,\lambda_i)=0$ implies $\lambda_i=0$. We now consider two
options,  regarding the dimension of the corresponding block
$C_{n_i}(\lambda_i)=C_{n_i}(0)$.\medskip

 \noindent{\bf Step 4}. Suppose $\lambda_i=0$ and the corresponding
block~$C_{n_i}(\lambda_i)$ has dimension~$n_i\ge 2$. For simplicity,
assume~$i=1$, so~$\lambda_1=0$.

Now, the first block of~$\mathfrak{p}(X,A,\dots,A)$ equals
$$X_{11}C_{n_1}(0)^{k-1}+
\bigl(\beta_2\,C_{n_1}(0)X_{11}C_{n_1}(0)^{k-2}+\dots+
 \beta_{k}\,C_{n_1}(0)^{k-1}X_{11} \bigr)$$
 Note that $C_{n_1}(0)$ is an upper-triangular
nilpotent. Since~$n_1\ge 2$, a straightforward computations show that
$E_{1n_1}C_{n_1}(0)^{k-1}+
\bigl(\beta_2\,C_{n_1}(0)E_{1n_1}C_{n_1}(0)^{k-1}+\dots+
 \beta_{k}\,C_{n_1}(0)^{k-1}E_{1n_1} \bigr)=0$.  Therefore,
$X:=E_{1n_1}\in\Omega_{\bullet A}$. {Similarly, we can also show that
$E_{1n_1}\in\Omega_{A\bullet A}$} and we are in the condition (ii)
again.\medskip

\noindent {\bf Step 5}. Finally, suppose  we {\em are not} under the
conditions of Step 2 and there exists $i$ such that~$\lambda_i=0$ with
the corresponding block~$C_{n_i}(\lambda_i)=C_{n_i}(0)$ of
dimension~$n_i=1$. Again, for simplicity $i=1$, so that $A=0\oplus C$,
where~$C:= C_{n_2}(\lambda_2)\oplus\dots\oplus C_{n_r}(\lambda_r)$ is
an $(n-1) \times (n-1)$ matrix. Recall
that~$\lambda_2,\dots,\lambda_r\not=0$, {so~$C$ is invertible}.

It is straightforward to see that $\FF E_{11}\subseteq \Omega_{\bullet
A}\cap \Omega_{A\bullet A}$ in this case. Let us show that {also} $\FF
E_{11}\supseteq \Omega_{\bullet A}\cap \Omega_{A\bullet A}$.\medskip

Retaining the block structure of~$X$, the
Eqs.~(\ref{eq:pom4})--(\ref{eq:pom4.1}) simplify for the blocks in the
first row/column into $X_{1j} C_{n_j}(\lambda_j)^{k-1}=0$,
respectively, into $\xi_{\jmath_0k}\,C_{n_j}(\lambda_j)^{k-1}X_{j1}
=0$. Since $C_{n_j}(\lambda_j)$ are invertible for $j=2,\dots,r$, and
since $\xi_{\jmath_0k}\not=0$, we get $X_{1j}=0=X_{j1}$ whenever $j\ge
2$.

Consider also the block $X_{st}$ for $s,t\ge 2$. Now, if $s\not=t$, it
is impossible to have simultaneously $p_{\bullet
A}(\lambda_s,\lambda_t)=0=p_{ A\bullet A}(\lambda_s,\lambda_t)$, by
conditions of Step~2.  This remains true if $s=t\ge 2$, for otherwise
$p_{\bullet A}(\lambda_s,\lambda_s)=0$, which would wrongly imply
$\lambda_s=0$. Then, however, we may copy the arguments from Step~1 to
deduce $X_{st}=0$. Therefore, the only possible nonzero block of~$X$ is
$X_{11}$, and so $X=\alpha E_{11}$ for some scalar~$\alpha$. Therefore,
$\Omega_{\bullet A}\cap\Omega_{A\bullet A }=\FF E_{11}$. This gives the
case~(iii) with $P=E_{11}$.
\end{proof}

Recall that ${\cal I}^1$ is the set of rank-one idempotents in ${\cal
M}_n(\FF)$.

\begin{corollary}\label{cor:Phi-perserves-rk1-idempote-in-both-directions}%
Let conditions of Theorem~\ref{thm:main-nonvanishing-sum} be satisfied.
Then
\begin{itemize}
\item[(i)] $\Phi({\cal I}^1)\subseteq \FF {\cal I}^1$.
\item[(ii)] If $\Phi(X)\in {\cal I}^1$ then $X=\alpha Y$ for certain
$\alpha\in \FF\setminus\{0\}$ and $Y\in {\cal I}^1$.
\end{itemize}
\end{corollary}

\begin{proof}
{\bf Step 1}. Pick a rank-one idempotent~$P$.  Then, $A:=\Id-P$ is an
idempotent of rank $(n-1)$. Thus, $A\in\DPhi_1$. The direct
calculations show that $\Omega_{\bullet A}\cap \Omega_{A\bullet
A}=\Omega_{\bullet(\Id-P)}\cap \Omega_{(\Id-P)\bullet (\Id-P)}$
consists precisely of those matrices~$X$ which satisfy
$$XA+\beta AXA+\beta_{k}
AX=0=\xi_{\jmath_0k}\,AX+\tilde{\beta} AXA+\tilde{\beta}_{1} XA;
\quad(\beta:=\beta_2+\dots+\beta_{k-1}).$$ By
Lemma~\ref{lem:A-circ-B=0<=>AB=0=BA}, $X$ is orthogonal to
idempotent~$A=\Id-P$. Hence, $X=\lambda P$, and so, $\Omega_{\bullet
A}\cap \Omega_{A\bullet A}=\FF\,P$.\medskip

{\bf Step 2}. Since $\Phi$ preserves zeros of $\mathfrak{p}$, the first
step implies $\Phi(P)\in\Omega_{\bullet\Phi(A)}\cap
\Omega_{\Phi(A)\bullet \Phi(A)}$ . By
Corollary~\ref{cor:Phi-preserves-0},~$\Phi(P)\not=0$, so
$\Omega_{\bullet\Phi(A)}\cap \Omega_{\Phi(A)\bullet
\Phi(A)}\not=\{0\}$. By Lemma~\ref{lem:charact_of_minimal_idempotents}
it follows that $\Omega_{\bullet\Phi(A)}\cap \Omega_{\Phi(A)\bullet
\Phi(A)}$ either is equal to a scalar multiple of an idempotent or
contains a rank-one nilpotent. We assume erroneously  the later, i.e.,
that there exists a rank-one nilpotent $Y\in\Omega_{\bullet\Phi(A)}\cap
\Omega_{\Phi(A)\bullet \Phi(A)}$. By the surjectivity of $\Phi$ it
follows that $Y=\Phi(X)$ for some~$X\in \DPhi_1$. Since $\Phi$
preserves zeros of $\mathfrak{p}$ strongly  we have
$X\in(\Omega_{\bullet A}\cap\Omega_{A\bullet A })\cap\DPhi_1=\FF P\cap
\DPhi_1$.

By Corollary~\ref{cor:Phi-preserves-0}, $X\not=0$. On the other hand,
$Y^2=0$ since $Y$ is a nilpotent of rank-one. {Hence},
$\mathfrak{p}(Y,\dots, Y)=0$, so also
$0=\mathfrak{p}(X,\dots,X)=(1+\xi)X^k$. This is clearly a contradiction
since $X\in\FF P$. Hence, $\Omega_{\bullet\Phi(A)}\cap
\Omega_{\Phi(A)\bullet \Phi(A)}$ contains no rank-one nilpotents.
{Therefore}, $\Phi(P)\in\Omega_{\bullet\Phi(A)}\cap
\Omega_{\Phi(A)\bullet \Phi(A)}=\FF Q$, for some rank-one
idempotent~$Q$. This proves~(i).\medskip

{\bf Step 3}. Conversely, if~$\Phi(X)$ is an idempotent of rank-one,
then $\FF\,\Phi(X)=\Omega_{\bullet B}\cap\Omega_{B\bullet B}$ for
$B:=\Id-\Phi(X)\in\DPhi_2$. Since $\Phi$ is surjective and preserves
zeros of $\mathfrak{p}$ strongly, we can prove that $X\in\FF\, P$ for
some rank-one idempotent~$P$  in a similar way, as in the proof of Item
{(i)}.
\end{proof}

\begin{remark}\label{rem_redefined_Phi}%
Since the polynomial $\mathfrak{p}$ is homogeneous, the assumptions and
the conclusion of Theorem~\ref{thm:main-nonvanishing-sum} will not be
affected if we replace~$\Phi$ by a mapping $\hat \Phi:A\mapsto
\delta(A)\cdot\Phi(A)$, where~$\delta: {\DPhi_1}\to \FF\backslash\{0\}$
is a scalar function. We may define~$\delta$ in such a way
that~$\hat{\Phi}$ preserves rank-one idempotents (i.e., $\hat
\Phi({\cal I}^1)\subseteq {\cal I}^1$), while~$\hat{\Phi}(A)=\Phi(A)$
for any other matrix from~$\DPhi_1$.

The redefined~$\hat{\Phi}$ may not be surjective, but we clearly have
${\cal I}^1\subseteq\im\hat{\Phi}$.
\end{remark}\bigskip

We are ready now to prove the main result.
\begin{proof}[Proof of Theorem~\ref{thm:main-nonvanishing-sum}]
We proceed in several steps.\medskip

{\bf Step 1}. {\em The transformation~$\hat\Phi$ preserves
orthogonality among rank-one idempotents.\/}

Indeed, suppose~$P,Q$ are two orthogonal rank-one idempotents. Then
$Q\in\Omega_{\bullet P}\cap\Omega_{P\bullet P}$. Therefore,
$\hat\Phi(Q)\in\Omega_{\bullet\hat\Phi(P)}\cap
\Omega_{\hat\Phi(P)\bullet \hat\Phi(P)}$, so that $X\hat\Phi(P)+\beta
\hat\Phi(P)X\hat\Phi(P)+\beta_{k} \hat\Phi(P)X=0
=\hat\Phi(P)X+\hat{\beta} \hat\Phi(P)X\hat\Phi(P)+\hat{\beta}_{1}
X\hat\Phi(P)$ for $X:=\hat\Phi(Q)$. The orthogonality between
$X=\hat\Phi(Q)$ and $\hat\Phi(P)$ now follows from
Lemma~\ref{lem:A-circ-B=0<=>AB=0=BA}.\medskip

{\bf Step 2.} {\em Let us show that $\hat\Phi$ is injective
transformation on the set of rank-one idempotents.\/}

We assume erroneously that~$\hat{\Phi}(P_1)=\hat{\Phi}(P_2)$ for some
{\em distinct} rank-one idempotents $P_1={\bf x}_1{\bf f}_1^{\tr}$ and
$P_2={\bf x}_2{\bf f}_2^{\tr}$. Then either ${\bf x}_1,{\bf x}_2$ are
linearly independent or ${\bf f}_1,{\bf f}_2$ are linearly independent
or both. Assume that ${\bf x}_1,{\bf x}_2$ are. Then, we can construct
a rank-one idempotent~$Q$ such that
$\mathfrak{p}(Q,P_1,\dots,P_1)=QP_1+\beta P_1QP_1+\beta_k P_1Q=0$,
but~$\mathfrak{p}(Q,P_2,\dots,P_2)=QP_2+\beta P_2QP_2+\beta_k
P_2Q\not=0$.

Indeed, we choose any nonzero ${\bf y}$ with
 ${\bf f}_1^\tr{\bf y}=0={\bf f}_2^\tr{\bf y}$.
 Suppose first  ${\bf y}=\mu_1 {\bf x}_1+\mu_2 {\bf x}_2$.
Then, $\mu_2\not=0$, since otherwise ${\bf y}={\mu_1} {\bf x}_1$,
contradicting ${\bf f}_1^{\tr}{\bf x}_1=1$, ${\bf f}_1^{\tr}{\bf y}=0$.
Now, as ${\bf x}_1,{\bf x}_2$ are linearly independent, we may choose
${\bf g}$ such that $ {\bf g}^\tr{\bf x}_1=0$ and ${\bf g}^\tr{\bf x}_2
=1/\mu_2$. Then  ${\bf g}^\tr{\bf y}=1$. Now, if ${\bf y}, {\bf
x}_1,{\bf x}_2$ are linearly independent, we may choose ${\bf g}$ such
that ${\bf g}^\tr{\bf x}_1 =0$, ${\bf g}^\tr{\bf x}_2  =1$, ${\bf
g}^\tr{\bf y}=1$. In both cases, $Q:={\bf
y} {\bf g}^{\tr}$ is the required idempotent. \\
For the chosen $Q$ we have
$$0=\mathfrak{p}(\hat\Phi(Q),\hat\Phi(P_1),\dots,\hat\Phi
(P_1))=\mathfrak{p}(\hat\Phi(Q),\hat\Phi(P_2),\dots,\hat\Phi
(P_2))\not=0,$$ a contradiction.

On the other hand, if~${\bf f}_1,{\bf f}_2$ are independent, we can
similarly find~$Q$ with $\mathfrak{p}(P_1,Q,\dots,Q)=0$,
but~$\mathfrak{p}(P_2,Q,\dots,Q)\not=0$. As before, this leads to a
contradiction. {Indeed}, $\hat\Phi(P_1)\ne \hat\Phi(P_2)$.
\medskip

{\bf Step 3.} {\em Now we can characterize the action of $\hat\Phi$ on
the set of rank-one idempotents.\/}

Injective mappings on rank-one idempotents, which preserve
orthogonality are classified in
\v{S}emrl~\cite[Theorem~2.3]{semrl:commutativity} ({\small this result
is stated only for~$\FF=\mathbb{C}$, but it was already remarked by the
author that the proofs are valid for {any algebraically closed field}
}). It follows that there exists a field
homomorphism~$\varphi:\FF\to\FF$, and invertible matrix~$T$ such that
either $\hat\Phi(P)=TP^\varphi T^{-1}$ for every rank-one
idempotent~$P$, or else $\hat\Phi(P)=T(P^\varphi)^{\tr} T^{-1}$ for
every rank-one idempotent~$P$.\medskip

{\bf Step 4.} {\em Field homomorphism $\varphi$ is surjective.\/}

It suffices to see that the restriction of $\hat\Phi$ on the set of
rank-one idempotents,  $\hat\Phi|_{{\cal I}^{1}}:{\cal I}^1\to{\cal
I}^1$, is surjective. Let us take any ${F}\in {\cal I}^1$. By
Remark~\ref{rem_redefined_Phi}, ${\cal I}^1\subseteq\im\hat\Phi$, thus
${F}=\hat\Phi(X)$ for certain $X\in \DPhi_1$.
Corollary~\ref{cor:Phi-perserves-rk1-idempote-in-both-directions} shows
that $X=\lambda P$ for some $P\in{\cal I}^1$. Now, choose pairwise
orthogonal $P_2,\dots,P_n\in{\cal I}^1$ that are also orthogonal
to~$P$. Clearly, $(\lambda P)\in\Omega_{\bullet P_i}\cap
\Omega_{P_i\bullet P_i}$, so also $\hat\Phi(\lambda
P)\in\Omega_{\bullet \hat\Phi(P_i)}\cap \Omega_{\hat\Phi(P_i)\bullet
\hat\Phi(P_i)}$. As in Step~1 we derive that ${F}=\hat\Phi(\lambda P)$
is orthogonal to~$\hat\Phi(P_2),\dots,\hat\Phi(P_n)$. On the other
hand, however, $\hat\Phi(P),\hat\Phi(P_2),\dots,\hat\Phi(P_n)$ are~$n$
pairwise orthogonal rank-one idempotents as well. This is possible only
when $\hat\Phi(P)={F}$, and the result follows.
\medskip

{\bf Step 5.} {\em The conclusion of
Theorem~\ref{thm:main-nonvanishing-sum} is valid for all {\em
non-nilpotent} rank-one matrices.\/}

Similar to Step 4 it can be shown that~$\hat\Phi(\lambda P) \in \FF\,
\hat\Phi(P)$ for any $P\in{\cal I}^1$, i.e., there exists a
transformation $\delta' : {\DPhi_1 \to \FF\backslash\{0\}}$ such that
$\hat\Phi(\lambda P)=\delta'(\lambda P)\hat\Phi(P)$ for all
$\lambda\in\FF$, $P\in {\cal I}^1$.
\medskip

{\bf Step 6.} {\em $\hat\Phi$ preserves the set of rank-one
 nilpotents.\/}

 To see this, we choose any rank-one nilpotent~$N$. Using similarity, we
may assume~$N=E_{12}$. Then, we can find $n-2$ pairwise orthogonal
rank-one idempotents $E_{33},\dots,E_{nn}$ that are also orthogonal
to~$N$. Clearly, $N\in\Omega_{\bullet E_{ii}}\cap \Omega_{E_{ii}\bullet
E_{ii}}$, so also $\hat\Phi(N)\in\Omega_{\bullet \hat\Phi(E_{ii})}\cap
\Omega_{\hat\Phi(E_{ii})\bullet \hat\Phi(E_{ii})}$ for $n-2$ pairwise
orthogonal rank-one idempotents~$\hat\Phi(E_{ii})$, $i=3,\dots,n$.

Using similarity in the image space, we may
assume~$\hat\Phi(E_{ii})=E_{ii}$. Thus, we have
$\hat\Phi(N)\in\Omega_{\bullet E_{ii}}\cap \Omega_{E_{ii}\bullet
E_{ii}}$ for $i=3,\dots,n$. As in Step~1 we derive that
idempotents~$E_{ii}$ are orthogonal on $\hat\Phi(N)$ for $i=3,\dots,n$.
Consequently,~$\hat\Phi(N)$ could be nonzero only at the upper left $2
\times 2$ block. On the other hand,
$\mathfrak{p}(N,\dots,N)=(1+\xi)N^k=0$. Thus
also~$0=1/(1+\xi)\mathfrak{p}(\hat\Phi(N),\dots,\hat\Phi(N))$, i.e.,
$\hat\Phi(N)^k=0$. Hence,~$\hat\Phi(N)$ is a nonzero nilpotent and all
its non-zero elements are concentrated in the $2\times 2$ upper-left
block. Thus $\hat\Phi(N)$ is a nilpotent of rank-one.\medskip

{\bf Step 7.} {\em The end of the proof}.

 Finally, consider the redefined
$\tilde{\Phi}:A\mapsto T^{-1}\hat\Phi (A)T$. By Step~3 either
$\tilde{\Phi}(P)\equiv P^\varphi$ for rank-one idempotents~$P$, or else
$\tilde{\Phi}(P)\equiv (P^\varphi)^{\tr}$ for rank-one idempotents~$P$.
Hence, applying~Lemma~\ref{lem:charact_(nilpot-of-rk1)_with_idempot}
({\small respectively,
Lemma~\ref{lem:charact_(nilpot-of-rk1)_with_idempot'}}) to rank-one
nilpotents $N_1:=N$ and $N_2:=\tilde{\Phi}(N)$ we obtain
$\tilde{\Phi}(N)=\delta''(N)\cdot N^\varphi$ for certain $\delta'':
{\cal M}_n(\FF)\to \FF$ ({\small respectively,
$\tilde{\Phi}(N)=\delta''(N)\cdot (N^\varphi)^{\tr}$}). Obviously, this
holds for any rank-one nilpotent~$N$.
\end{proof}

\subsection{Proof of Corollaries}

It will be beneficial to regard $\FF^n$ as the space of  matrices of
dimension $n$--by--$1$, i.e., {\em column vectors}. Thus, any rank--one
matrix~$A\in{\cal M}_n(\FF)$ can be written as  $A={\bf x}{\bf
f}^{\tr}$ for suitably chosen  ${\bf x},{\bf f}\in\FF^n$. Its trace
then equals $\mathrm{Tr}\, A={\bf f}^{\tr}{\bf x}$.

Now, to prove Corollary~\ref{cor:TAllMatr}, we will rely on the
folowing result due to Bre\v{s}ar and
\v{S}emrl~\cite[Thm.~2.4]{bresar_semrl}, which we state slightly
changed, recasting it into our framework:
\begin{lemma}\label{lem:bresar_semrl}%
Let  $\FF$ be an infinite field  with $\ch\FF\not= 2$, and let
$R_1,R_2,R_3 \in{\cal M}_n(\FF)$ be three matrices. Then (i) implies
(ii) below.
\begin{itemize}
\item[(i)]  The vectors $R_1{\bf u}$, $R_2{\bf
u}$, and $R_3{\bf u}$ are linearly dependent for every ${\bf u}\in
\FF^n$.
\item[(ii)]Either $R_1,R_2,R_3$ are linearly dependent, or there
exist ${\bf v},{\bf w},{\bf z} \in \FF^n$ such that $\Im R_i\subseteq
\Lin_{\FF}\{{\bf v},{\bf w},{\bf z}\}$ for $i = 1,2,3$, or there exists
a rank-one idempotent $P\in{\cal M}_n(\FF)$ such that
$$\dim\Lin\nolimits_{\FF}\{(\Id - P)R_1, \;(\Id - P)R_2,\; (\Id - P)R_3\} = 1.$$
\end{itemize}
\end{lemma}

With its help, the following generalization of
Lemma~\ref{lem:charact_(nilpot-of-rk1)_with_idempot} can be proven:

\begin{lemma}\label{lem:A=lambda*B_iff...}%
Let~$n\ge 4$ be an integer, let $A,B\in {\cal M}_n(\FF)$ be nonzero
matrices, let $\varphi:\FF\to\FF$ be a nonzero field homomorphism, and
let~$\alpha,\beta\in\FF$. Assume that the following condition (i) is
satisfied:
\begin{itemize}
\item[(i)] {$NAP+\alpha PAN=0\Longleftrightarrow
    0=N^\varphi B P^\varphi+\beta P^\varphi B N^\varphi$  holds for every
    rank-one idempotent~$P$ and every rank-one matrix~$N$
    with~$PN=0=NP$.}
\end{itemize}
Then there exists $\gamma,\mu\in\FF$ such that $B=\gamma
A^\varphi+\mu\Id$.
\end{lemma}
\begin{proof}
Pick any nonzero vector  ${\bf x}\in\FF^n$ and assume erroneously that
the vector ${\bf b}:=B{\bf x}^\varphi\not\in\Lin_{\FF}\{A^\varphi {\bf
x}^\varphi,{\bf x}^\varphi\}$. Denote ${\bf a}:=A{\bf x}$, and let
${\bf f}_1,\dots,{\bf f}_{l}$ be a basis of~$\{{\bf a},{\bf
x}\}^\bot:=\{{\bf f}\in\FF^n;\;\; {\bf f}^{\tr}{\bf a}=0={\bf
f}^{\tr}{\bf x}\}$ ({\small here,~$l=n-2$ or $n-1$  if ${\bf x}$ and
${\bf a}$ are linear independent or not, correspondingly}). Since the
rank of a matrix equals the maximal size of its nonzero minors, the
vectors~${\bf f}_1^\varphi,\dots,{\bf f}_{l}^\varphi$ are also linearly
independent. Hence, they form a basis of~$\{{\bf a}^\varphi,{\bf
x}^\varphi\}^\bot$. Now,~${\bf b}\not\in\Lin_{\FF}\{{\bf
a}^\varphi,{\bf x}^\varphi\}$, so~$({\bf f}_j^\varphi)^{\tr}{\bf
b}\not=0$ for at least one~$j$. Consequently, there exists ${\bf
f}={\bf f}_j$ such that
$${\bf f}^{\tr}{\bf x}=0={\bf f}^{\tr}A{\bf x}, \quad\hbox{and}\quad ({\bf
f}^\varphi)^{\tr}B{\bf x}^\varphi\not=0.$$
 Since~$n > 2$ we can find~${\bf y}$
such that
\begin{equation} \label{Eq:x_y}
{\bf x}\not\in\Lin\nolimits_{\FF}\{{\bf y},A{\bf y}\}.
\end{equation}
 Indeed, write
$\FF^n=\Lin_{\FF}\{{\bf x}\}\oplus M$. If $\Ker (A|_M)\not=0$, then any
nonzero ${\bf y}\in\Ker (A|_M)$ satisfies Eq.~(\ref{Eq:x_y}). Assume
now that $\Ker (A|_M)=0$ and ${\bf x}\in\Lin_{\FF}\{{\bf y},A{\bf y}\}$
for each~${\bf y}\in M$. Then $A|_M{\bf y}=\lambda_{\bf y}{\bf
x}+\mu_{\bf y}{\bf y}$. Since ${\bf x}\notin M$, we could deduce that
$\lambda_{\bf y}$ is a linear functional, on the space~$M$ with $\dim
M\ge 2$. Hence, $\lambda_{\bf y}=0$ for at least one nonzero ${\bf
y}={\bf y}_0\in M$. For this ${\bf y}_0$ we have $A{\bf y}_0=\mu_{{\bf
y}_0}{\bf y}_0$ and then $\Lin_{\FF}\{{\bf y}_0,A{\bf y}_0\}=\FF {\bf
y}_0$. However, ${\bf x}\notin\FF {\bf y}_0 \subset M$ --- a
contradiction.\medskip

Now, by~(\ref{Eq:x_y}), we may choose a vector~$\bf g$ with ${\bf
g}^{\tr}{\bf y}=0={\bf g}^{\tr} A{\bf y}$, but ${\bf g}^{\tr}{\bf
x}=1$. Then, $P:={\bf x}{\bf g}^{\tr}$ is an idempotent of rank-one,
and $N:={\bf y}{\bf f}^{\tr}$ is a matrix of rank-one, and we have
$PN=0=NP$. Moreover, $NAP+\alpha PAN= ({\bf f}^{\tr}A{\bf x})\,{\bf
y}{\bf g}^{\tr}+\alpha ({\bf g}^{\tr}A{\bf y})\,{\bf x}{\bf f}^{\tr}=
0+\alpha\cdot 0=0$. Consequently,  the condition~(i) implies
\begin{align*}
0&= ( ({\bf f}^\varphi)^{\tr}B{\bf x}^\varphi)\cdot {\bf y}^\varphi
({\bf g}^\varphi)^{\tr}+
 \beta\cdot (({\bf g}^\varphi)^{\tr} B{\bf y}^\varphi)\cdot {\bf x}^\varphi({\bf
 f}^\varphi)^{\tr}.
\end{align*}
However, the first summand on the right is nonzero. Hence, the right
side is nonzero, since ${\bf y}^\varphi({\bf g}^\varphi)^{\tr} $ and
${\bf x}^\varphi({\bf f}^\varphi)^{\tr}$ are linearly independent
matrices ({\small namely, $({\bf g}^\varphi)^{\tr}{\bf x}^\varphi=1$,
while $({\bf f}^\varphi)^{\tr}{\bf x}^\varphi=0$}). This contradiction
establishes that
\begin{equation}\label{eq:Bx_in_(x,Ax)}%
B{\bf x}^\varphi\in\Lin\nolimits_{\FF}\{A^\varphi{\bf x}^\varphi,{\bf
x}^\varphi\}\qquad \hbox{ for any ${\bf x}\in\FF^n$.}
\end{equation} \medskip

By Eq.~(\ref{eq:Bx_in_(x,Ax)}), the vectors $B{\bf x}^\varphi, {\bf
x}^\varphi, A^\varphi {\bf x}^\varphi$ are always $\FF$--linearly
dependent. Let us show that even more is true: indeed the matrices
$B,\Id,A^\varphi$ are locally linearly dependent, i.e., for any ${\bf
z}\in \FF^n$ the vectors $B{\bf z},{\bf z}, A^\varphi {\bf z}$ are
linearly dependent. To demonstrate this, we consider a matrix $\Xi({\bf
z}):=[B{\bf z},{\bf z},A^\varphi{\bf z}]$ with three columns. By
Eq.~(\ref{eq:Bx_in_(x,Ax)}), if ${\bf z}={\bf x}^\varphi$ for a certain
${\bf x}\in\FF^n$, then all its $3$--by--$3$ minors must vanish.

Consider any such minor. It is a polynomial $q(z_1,\ldots, z_n)\in
\FF[z_1,\ldots, z_n]$, where ${\bf z}=[z_1,\ldots, z_n]^{\rm tr}$. By
Eq.~(\ref{eq:Bx_in_(x,Ax)}) this polynomial vanishes identically
whenever the variables take the values from a subfield ${\cal
O}:=\varphi(\FF)\subseteq\FF$, i.e., for any values $\alpha_1, \ldots,
\alpha_n \in {\cal O}$ it holds that $q(\alpha_1,\ldots, \alpha_n)=0$.
Now, being algebraically closed,~$\FF$ and hence also ${\cal
O}=\varphi(\FF)$ have infinitely many  elements. It is easy to see then
that then,~$q$ is a zero polynomial. {\small For the sake of
completeness we  sketch  the proof here. We will write $q$ in the form
$$q(z_1,\dots,z_n)=a_n(z_1,\dots,z_{n-1})z_n^n
+\dots+a_1(z_1,\dots,z_{n-1})z_n+ a_0(z_1,\dots,z_{n-1}).$$ By the
assumptions, this vanishes whenever $z_1,\dots,z_n\in{\cal O}$. Now, at
each fixed $z_1,\dots,z_{n-1}\in{\cal O}$, this is a polynomial in only
one variable, $z_n$. However, it is zero for infinitely many values of
$z_n\in{\cal O}$. Hence, $x\mapsto q(z_1,\dots,z_{n-1},x)$ is a zero
polynomial for each fixed $(z_1,\dots,z_{n-1})\in {\cal O}^{n-1}$. That
is, all its coefficients, $a_i(z_1,\dots,z_{n-1})$  are zero for any
$(z_1,\dots,z_{n-1})\in{\cal O}^{n-1}$. It remains to show that
$a_i(z_1,\dots,z_{n-1})$ are identically zero, not only for any choice
of $z_1,\dots,z_{n-1}\in {\cal O}$, but also for any choice of
$z_1,\dots,z_{n-1}\in\FF$. Now,
we may repeat the aforesaid procedure with each $a_i(z_1,\dots,z_{n-1})$:
Write it as
$$a_i(z_1,\dots,z_{n-1})=b_{im}(z_1,\dots,z_{n-2})z_{n-1}^m+\dots  +
    b_{i0}(z_1,\dots,z_{n-2})$$
and argue as before to deduce that $b_j(z_1,\dots,z_{n-2})$ vanishes
for any choice of $z_1,\dots,z_{n-2}\in{\cal O}$. Continuing in this
way we obtain at the end certain polynomials $c_l(z_1)\in\FF[z_1]$
which are zero for any value $z_1\in {\cal O}$. It follows that
$c_l(z_1)$ is zero for infinitely many values of $z_1$, i.e., that
$c_l(z_1)$ is a zero polynomial. By the backward induction, we get that
all coefficients~$a_i(z_1,\dots,z_{n-1})$ are  zero, i.e., that $q$ is
indeed a zero polynomial.    }

Therefore, $q(z_1,\dots,z_n)=0$ holds for any $z_1,\ldots,z_n\in\FF$.
We repeat this with all $3$--by--$3$ minors to deduce that $\rk\Xi({\bf
z})\le 2$ for any ${\bf z}\in\FF^n$, as claimed.\medskip

Consequently,  $(R_1,R_2,R_3):=(B,\Id,A^\varphi)$ are locally linearly
dependent. We can now invoke Bre\v{s}ar and \v{S}emrl's theorem, see
Lemma~\ref{lem:bresar_semrl} in this text. Since~$R_2=\Id$ and
$\dim(\im(\Id))=n\gneq 3$, the only three possibilities left to
consider are (a) $A^\varphi=\lambda\Id$, or (b)
$B=\gamma\,A^\varphi+\mu\Id$, or
\begin{equation}
(\Id-Q)B=\lambda(\Id-Q)\Id,\qquad
(\Id-Q)A^\varphi=\lambda'(\Id-Q)\Id\tag{c}
\end{equation} for some rank-one idempotent~$Q$.
Under (a), $B$ must also be a scalar, in view of
Eq.~(\ref{eq:Bx_in_(x,Ax)}). So, under (a)--(b) we are done.\medskip

Consider lastly~(c). Decomposing $B=(\Id-Q)B+QB$, and using Eq.~(c),
easily reveals $B=\lambda\Id+Q(B-\lambda\Id)=\lambda\Id+\hat{\bf
u}\hat{\bf v}_B^{\tr}$ for some column vectors $\hat{\bf u},\hat{\bf
v}_B$. Likewise, $A^\varphi=\lambda'\Id+\hat{\bf u}\hat{\bf
v}_A^{\tr}$. By passing the appropriate scalar to the other term (in
$\hat{\bf u}\hat{\bf v}_B^{\tr}$, $\hat{\bf u}\hat{\bf v}_A^{\tr}$), we
may assume that at least one entry of vector~$\hat{\bf u}$ equals~$1$.
Then, from $\lambda'\Id+\hat{\bf u}\hat{\bf
v}_A^{\tr}=A^\varphi\in{\cal M}_n(\varphi(\FF))$ it follows that
$\hat{\bf u},\hat{\bf v}_A\in\varphi(\FF^n)$, and also
$\lambda'\in\varphi(\FF)$. Let ${\bf u}$, ${\bf v}_A$, and $\lambda''$
be such that $\hat{\bf u}={\bf u}^{\varphi} $, $\hat{\bf v}_A={\bf
v}_A^{\varphi}$, and $\lambda'=\varphi(\lambda'')$. Then
$$A^\varphi=\varphi(\lambda'')\Id+{\bf u}^\varphi({\bf
v}_A^\varphi)^{\tr} \quad\hbox{and}\quad B=\lambda\Id+{\bf
u}^\varphi{\bf v}_B^{\tr}.$$
 Now, if ${\bf v}_A^\varphi,{\bf v}_B$ are linearly dependent we are
done. Assume erroneously they are not.  We first choose ${\bf v}$ such
that ${\bf v}^\varphi,{\bf v}_A^\varphi,{\bf v}_B$ are independent, and
at the same time, $({\bf v}^\varphi)^{\tr}{\bf u}^\varphi=1$ ({\small
Such a vector ${\bf v}$ exists since we can enlarge~${\bf v}_A$ with
${\bf v}_2,\dots,{\bf v}_{n} $ to a basis of~$\FF^n$, assuming ${\bf
v}_i^{\tr}{\bf u}=1$. Then, ${\bf v}_A^\varphi$, ${\bf
v}_2^\varphi,\dots,{\bf v}_n^\varphi$ is still a basis, so some ${\bf
v}^\varphi:={\bf v}_i^\varphi$ is independent of~${\bf
v}_A^\varphi,{\bf v}_B$}).

By the choice of ${\bf v}$, the vector ${\bf u}^\varphi\notin \{{\bf
v}^\varphi\}^\bot$. So much the more ${\bf u}^\varphi\notin \{{\bf
v}^\varphi,{\bf v}_A^\varphi\}^\bot$, so that $\Lin_\FF\{{\bf
u}^\varphi\}\cap \{{\bf v}^\varphi,{\bf v}_A^\varphi\}^\bot=\{0\}$. We
next choose ${\bf w}$ with  ${\bf w}^\varphi\in\{{\bf v}^\varphi,{\bf
v}_A^\varphi\}^\bot\backslash\{{\bf v}_B\}^\bot $ ({\small it is
possible since ${\bf v}^\varphi,{\bf v}_A^\varphi,{\bf v}_B$ are
independent .}) Hence, ${\bf w},{\bf u}$ are independent, so we can
choose ${\bf h}$ with ${\bf h}^{\tr}{\bf u}=0$ and ${\bf h}^{\tr}{\bf
w}:=1$. Lastly, choose nonzero ${\bf s}\in \{{\bf h}\}^\bot$.

We now form $N_1:={\bf s}{\bf v}^{\tr}$ and $P_1:={\bf w}{\bf
h}^{\tr}$. By its choice, ${\bf w}\in\{{\bf v},{\bf v}_A\}^\bot$ so it
follows $P_1N_1=0=N_1P_1$, and $P_1^2=P_1$. Moreover,
$N_1AP_1=0=P_1AN_1$. By assumptive condition~(i), we would have to have
$0=N_1^\varphi BP_1^\varphi+\beta P_1^\varphi B N_1^\varphi$. However,
it equals
$$
N_1^\varphi BP_1^\varphi+\beta P_1^\varphi B
N_1^\varphi=\underbrace{\bigl( ({\bf v}^\varphi)^{\tr}{\bf u}^\varphi
\bigr)}_{=1}\cdot\underbrace{({\bf v}_B^{\tr}{\bf
w}^\varphi)}_{\not=0}\cdot {\bf s}^\varphi({\bf
h}^\varphi)^{\tr}+\beta\cdot 0 \not=0.$$ This contradiction finally
establishes linear dependence between ${\bf v}_A^\varphi$ and ${\bf
v}_B$.
\end{proof}\bigskip

\begin{proof}[Proof of Corollary~\ref{cor:TAllMatr}]
Obviously,~$\Phi$ satisfies the hypothesis of
Theorem~\ref{thm:main-nonvanishing-sum}, hence also its conclusion.
That is, (i) and (ii) hold for rank-one matrices. Now, we may
replace~$\Phi$ by a mapping $X\mapsto \frac{1}{ \gamma(X)}\,
T^{-1}\Phi(X)T$ to achieve that either~$\Phi(X)\equiv X^\varphi$ holds
for all rank-one matrices, or else $\Phi(X)\equiv (X^\varphi)^{\tr}$
holds for all rank-one matrices. It remains to see that, modulo scalar
multiplication and scalar addition,  same holds for~$A\in\DPhi_1$ of
rank~$\ge 2$.\medskip

Assume first $\Phi(X)\equiv X^\varphi$ for all~$X$ of rank-one, and
let~$A\in\DPhi_1$ be of rank $\ge 2$. Now, by definition, our
polynomials satisfy $\sum_{\sigma(1)=1}\alpha_\sigma=1$. Note that
$$1=\sum_{\sigma(1)=1}\alpha_\sigma=\sum_{\sigma(1)=1\atop\sigma(2)=2}
\alpha_\sigma+\sum_{\sigma(1)=1\atop\sigma(3)=2}\alpha_\sigma+\dots+\sum_{\sigma(1)=
1\atop\sigma(k)=2}\alpha_\sigma,$$ so at least one summand on the right
is nonzero. Say,
$\tau:=\sum_{\sigma(1)=1\atop\sigma(j'_0)=2}\alpha_\sigma\not=0$.
Having found~$j'_0$, we next pick any rank-one idempotent~$P$, and any
rank-one matrix~$N$ with~$PN=0=NP$. Consider now
$\mathfrak{p}(N,P,\dots,P, \fbox{$A$}_{j'_0},P,\dots,P)$ with matrix
$A$ at the position $j'_0$. An easy argument reveals that
\begin{align*}
\mathfrak{p}(N,P,\dots,P,
\fbox{$A$}_{j'_0},P,\dots,P)&=\sum_{\sigma(1)=1\atop\sigma(j'_0)=2}\alpha_\sigma
NAP+
\sum_{\sigma(1)=k\atop\sigma(j'_0)=k-1}\alpha_\sigma PAN+\\
&\hphantom{\sum_{\sigma(1)=1\atop\sigma(j'_0)=2}\alpha_\sigma
NAP}\mbox{}+
\sum_{\text{the rest permut.}}\alpha_\sigma\cdot 0\\[5mm]
&= \tau (NAP+\alpha PAN)
\end{align*}
where $\alpha:=1/\tau
\sum_{\sigma(1)=k\atop\sigma(j'_0)=k-1}\alpha_\sigma$.  Similarly, the
value of
$\mathfrak{p}(\Phi(N),\Phi(P),\dots,\linebreak[3]\Phi(P),\fbox{$\Phi(A)$}_{j'_0},\Phi(P),\dots,\Phi(P))=
\mathfrak{p}(N^\varphi,
P^\varphi,\dots,P^\varphi,\fbox{$\Phi(A)$}_{j'_0},P^\varphi,\dots,P^\varphi)$
further equals $\tau(N^\varphi\Phi(A)P^\varphi+\alpha
P^\varphi\Phi(A)N^\varphi)$. Consequently,

$$ NAP+\alpha PAN =0
\Longleftrightarrow 0=N^\varphi\Phi(A)P^\varphi+\alpha
P^\varphi\Phi(A)N^\varphi,$$ so Lemma~\ref{lem:A=lambda*B_iff...}
gives~$\Phi(A)=\gamma(A) A^\varphi+\mu(A)\Id$.\medskip

Assume lastly $\Phi(X)\equiv (X^\varphi)^{\tr}$ for all~$X$ of
rank-one. Pick rank-one~$N,P$, with~$P^2=P$ and~$PN=0=NP$. Also,  pick
any~$A\in\DPhi_1$ and let~$B:=\Phi(A)$. We deduce, {similarly} as
before, that
\begin{equation}\label{eq:A^tr}%
\begin{aligned}
& &0&=NAP+\alpha PAN\\
&&&=\tfrac{1}{\tau}\,\mathfrak{p}(N,P,\dots,P,\fbox{$A$}_{j'_0},P,\dots,P)\qquad
\Longleftrightarrow\\
\Longleftrightarrow& & 0 &=\tfrac{1}{\tau}\,\mathfrak{p}(\Phi(N),\Phi(P),
\dots,\Phi(P),\fbox{$\Phi(A)$}_{j'_0},\Phi(P),\dots,\Phi(P))\\
&&&=\tfrac{1}{\tau}\,\mathfrak{p}\bigl((N^\varphi)^{\tr},(P^\varphi)^{\tr},
\dots,(P^\varphi)^{\tr},B,(P^\varphi)^{\tr},\dots,(P^\varphi)^{\tr}\bigr)\\
&&&=(N^\varphi)^{\tr}B(P^\varphi)^{\tr}+\alpha
(P^\varphi)^{\tr}B(N^\varphi)^{\tr}
\end{aligned}
\end{equation}
Choose~$(N,A,P):=(E_{22},E_{12},E_{11})$. Then, $PN=0=NP$ and
$B:=\Phi(A)=(E_{12}^\varphi)^{\tr}=E_{21}$. On one hand, this gives
$NAP+\alpha PAN=\alpha\,E_{12}$, and,  on the other,
$(N^\varphi)^{\tr}B(P^\varphi)^{\tr}+\alpha
(P^\varphi)^{\tr}B(N^\varphi)^{\tr}=E_{21}$. Consequently, the
equivalence~(\ref{eq:A^tr}) reads $0=\alpha E_{12}\Longleftrightarrow
0=E_{21}$, and so $\alpha\not=0$.

We may now rewrite equivalence~(\ref{eq:A^tr}) into
$$NAP+\alpha PAN=0\Longleftrightarrow 0=\bigl(
N^\varphi B^{\tr}P^\varphi
\bigr)^{\tr}+\tfrac{1}{\alpha}\bigl(P^\varphi B^{\tr}N^\varphi
\bigr)^{\tr}.$$ The right side is further equivalent to $0= N^\varphi
B^{\tr}P^\varphi +\frac{1}{\alpha}P^\varphi B^{\tr}N^\varphi$. Hence,
Lemma~\ref{lem:A=lambda*B_iff...}
gives~$\Phi(A)^{\tr}=B^{\tr}=\gamma(A) A^\varphi+\mu(A)\Id$.
\end{proof}\bigskip
\begin{proof}[Proof of Corollary~\ref{cor:corolary_to_Jordan_product}]
Assumption $\ch\FF\not=2$ ensures $\mathfrak{p}(x,y):=xy+yx$ is a
polynomial with nonvanishing sum of coefficients. Since both
mappings~$X\mapsto X^{\varphi}$ and~$X\mapsto X^{\tr}$ preserve the
zeros of~$\mathfrak{p}(x,y)$, we may replace~$\Phi$ by a mapping
$X\mapsto\bigl( \frac{1}{ \gamma(X)}\, T^{-1}\Phi(X)T
\bigr)^{\varphi^{-1}}$,  respectively, by $X\mapsto \bigl(
\bigl(\frac{1}{ \gamma(X)}\,
T^{-1}\Phi(X)T\bigr)^{\tr}\bigr)^{\varphi^{-1}}$ to achieve that~$\Phi$
leaves fixed all rank-one matrices. It remains to see that
$\Phi(A)=\gamma(A) A$ holds also for diagonalizable  matrices~$A$ with
the spectrum $\{\lambda,-\lambda\}$. In view of
Corollary~\ref{cor:Phi-preserves-0} we may assume
further~$A\not=0$.\medskip

Using a similarity~$S$, we may write  $A=S\bigl( \lambda\Id_{n_1}\oplus
(-\lambda)\Id_{n_2}\bigr)S^{-1}$. It is easy to see
that $\Omega_A:=\{X\in{\cal M}_n(\FF);\;\;XA+\,AX=0\}=S\left[\begin{smallmatrix} 0 & K\\
L & 0
\end{smallmatrix}\right] S^{-1}$ where $K,L$ are arbitrary matrices of an
appropriate size. Now, since~$\Phi$ fixes rank-one matrices, we have,
by the defining equation~(\ref{eq:poly-with-nonvanishing-sum}),
$X\Phi(A)+\Phi(A)X=0$
at least for each rank-one $X\in S\left[\begin{smallmatrix} 0 & K\\
L & 0
\end{smallmatrix}\right] S^{-1}$. Obviously, the set~$\Omega_{\Phi(A)}$ of all
matrices~$X$ with $X\Phi(A)+\Phi(A)X=0$ is a vector subspace of~${\cal
M}_n(\FF)$, so actually
$$
\Omega_{\Phi(A)}\supseteq S\left[\begin{smallmatrix} 0 & K\\
L & 0
\end{smallmatrix}\right] S^{-1}.$$
We now write $\Phi(A)=S\left[\begin{smallmatrix} U & V\\
W & Z
\end{smallmatrix}\right]S^{-1}$, and solve the identity
$$\left[\begin{smallmatrix} U & V\\
W & Z
\end{smallmatrix}\right]\left[\begin{smallmatrix} 0 & K\\
L & 0
\end{smallmatrix}\right]+\left[\begin{smallmatrix} 0 & K\\
L & 0
\end{smallmatrix}\right]\left[\begin{smallmatrix} U & V\\
W & Z
\end{smallmatrix}\right]\equiv 0\qquad (\forall\,K,\,\forall\, L).$$
Straightforward calculations
give~$V=0=W$, and $V=\mu\Id_{n_1}$, $Z=-\mu\Id_{n_2}$. Thus,
$\Phi(A)=S\,\diag(\mu,-\mu)S^{-1}=\frac{\mu}{\lambda} \cdot A$.
\end{proof}\bigskip

\begin{proof}[Proof of Corollary~\ref{cor:nonsurjective}] We first
prove that~$\Phi$ maps rank-one idempotents into scalar multiples of
rank-one idempotents, and preserves their orthogonality.

Indeed, let $P_1,\ldots, P_n$ be the set of $n$ pairwise orthogonal
rank-one idempotents.  Clearly then,
$\mathfrak{p}(P_i,P_j,\ldots,P_j)=0=\mathfrak{p}(P_j,P_i,P_j,
\ldots,P_j)=\dots=\mathfrak{p}(P_j,\ldots,P_j,P_i)$ for all $i\ne j$.
This implies a similar set of equations on their $\Phi$--images
$A_i:=\Phi(P_i)$ and~$A_j:=\Phi(P_j)$. We write them explicitly:
$$\begin{aligned}
\sum_{\sigma(1)=1}\alpha_\sigma A_i A_j^{k-1}&+\sum_{\sigma(1)=2}
\alpha_\sigma A_jA_i A_j^{k-2}&+\dots +\sum_{\sigma(1)=k}\alpha_\sigma
A_j^{k-1}A_i&=0\\
\sum_{\sigma(2)=1}\alpha_\sigma A_i A_j^{k-1}&+\sum_{\sigma(2)=2}
\alpha_\sigma A_jA_i A_j^{k-2}&+\dots +\sum_{\sigma(2)=k}\alpha_\sigma
A_j^{k-1}A_i&=0\\
\multispan4{\dotfill}\\
 \sum_{\sigma(k)=1}\alpha_\sigma A_i
A_j^{k-1}&+\sum_{\sigma(k)=2}\alpha_\sigma
A_jA_iA_j^{k-2}&+\dots+\sum_{\sigma(k)=k}\alpha_\sigma A_j^{k-1}A_i&=0
\end{aligned}$$
These may be regarded as a system of~$k$ homogeneous linear equations
in `variables' $A_j^s A_iA_j^{k-1-s}$. By the assumptions, the
coefficient matrix is invertible, so the only solution is $A_j^s
A_iA_j^{k-1-s}=0$ for each~$s$. In particular,
$A_iA_j^{k-1}=A_j^{k-1}A_i=0$. That is,
\begin{equation}\label{eq:A_i*A_j=0}%
\Im(A_j^{k-1})\subseteq \Ker(A_i);\qquad (i\not= j).
\end{equation}
 Moreover, $\mathfrak{p} (P_i,\ldots, P_i)=(1+\xi)P_i^k=(1+\xi)P_i\ne 0$
implies that $0\ne\mathfrak{p}(A_i,\ldots,A_i)=(1+\xi)A_i^k$, for any
$i$.\medskip

We can now follow the arguments from~\cite[Lemma~2.2]{chan_li_sze} of
Chan, Li, and Sze: Firstly, we claim that $\rk(A_i)=1$ for any $i$.
Suppose on a contrary that, say, $\rk(A_1)\ge 2$. Then,
$\dim\Ker(A_1)<n-1$ and we deduce from~(\ref{eq:A_i*A_j=0}) that $\dim
(\Im(A_2^{k-1}) +\ldots +\Im (A_n^{k-1}))<n-1$. Hence, there exists $j$
such that
\begin{equation}\label{eq:sum}
\Im(A_j^{k-1})\subseteq \sum\limits_{l=2,\ldots,n \atop l\ne j} \Im
(A_l^{k-1}).
\end{equation}
Indeed, otherwise, by the induction, $\dim\bigl( \Im(A_2^{k-1}) +\ldots
+\Im (A_n^{k-1})\bigr)\ge n-1$, which is a contradiction. Again,
by~(\ref{eq:A_i*A_j=0}), the right hand side of Eq.~(\ref{eq:sum}) is
contained in $\Ker(A_j)$. Thus $A^k_j=0$ which contradicts
$A^k_j=\mathfrak{p}(A_j,\dots,A_j)\not=0$ by the assumption (ii).
Therefore, $\rk(A_i)=1$ for any~$i$.

Since $A_i^k\ne 0$  there exist $\lambda_i\in \FF\backslash\{0\}$ and
an idempotent matrix $Q_i$ of rank-one such that $A_i=\lambda_i Q_i$.
Lastly, it follows from $A_iA_j=\frac{1}{\lambda_j^{k-2}}
A_iA_j^{k-1}=0$ and $A_jA_i=\frac{1}{\lambda_i^{k-2}} A_jA_i^{k-1}=0$
that~$Q_i$ are pairwise orthogonal.\medskip

 Thus, $\Phi$
maps orthogonal idempotents of rank-one into scalar multiples of
orthogonal idempotents of rank-one. We now redefine~$\Phi$ as in
Remark~\ref{rem_redefined_Phi}. The rest
--- with the sole exception of Step~4 --- follows directly the proof of
Theorem~\ref{thm:main-nonvanishing-sum}.
\end{proof}

\section{\label{S3} Polynomials with vanishing sums of coefficients}

\begin{de}\label{def:admissible}
Let~$k\ge 2$. A subset of nonidentical permutations $\Xi\subseteq{\cal
S}_k$ is called an  {\em admissible subset\/} if the following two
conditions are satisfied:
\begin{itemize}
\item [(i)] There exists $t\in \{1,\ldots , k\}$ such that each
        $\sigma\in\Xi$ fixes the first $t-1$ elements, but
        $\sigma(t)\not=t$  ({\small note that $t<k$, otherwise,~$\sigma$ would have
        to be identical permutation}).
\item[(ii)] There exist integers $w,u,v\in \{1,\ldots , k\}$, $u<v$, such that
 $\sigma(w)=v$ and $\sigma(w+1)=u$ for each $\sigma\in \Xi$.
 \end{itemize}
\end{de}
Note that in particular, $\sigma(w+1)<\sigma(w)$ for each $\sigma\in
\Xi$.

\begin{example}
We give two examples of admissible sets which are the most visualizing
on one hand and which show that there are many admissible sets among
the subsets of ${\cal S}_k$ on the other hand.
\begin{itemize}
\item $\Xi:=\{\sigma\}$ is  admissible subset whenever~$\sigma$ is nonidentical.
\item An admissible subset is also the subset of all permutations from~${\cal
S}_k$ that swap $1$ and $2$ ({\small take $t:=1=:w$, $u:=1$, $v:=2$} in
Definition~\ref{def:admissible}).
\end{itemize}
It is not hard to see that the cardinality of an admissible set, with
$t=1$, i.e., which fixes no initial elements, is either $(k-2)!$ or
$(k-2)! - (k-3)!$, depending on choosing $w$ and $u,v$.
\end{example}

The main result of the present section is the following theorem. In
contrast to the Theorem~\ref{thm:main-nonvanishing-sum}  we do not
assume that $\FF$ is algebraically closed, and~$\Phi$ is a strong
preserver in this section.

\begin{theorem}\label{thm:main-on-vanishing-sum}%
Let $\FF$ be an arbitrary field with more than 2 elements, let $n\ge
3$,  $k\ge 2$ be integers, and let $\Xi\subset {\cal S}_k$ be a fixed
{\em admissible subset} of permutations. Suppose that
 two given homogeneous multilinear polynomials
$$\mathfrak{p}_1(x_1,\dots ,x_k):=x_1\cdots x_k-\sum_{\sigma\in\Xi} \alpha_\sigma
x_{\sigma(1)}\cdots x_{\sigma(k)};\qquad (\alpha_\sigma\in\FF)$$
$$\mathfrak{p}_2(x_1,\dots ,x_k):=x_1\cdots x_k-\sum_{\sigma\in\Xi} \beta_\sigma
x_{\sigma(1)}\cdots x_{\sigma(k)};\qquad (\beta_\sigma\in\FF)$$ satisfy
$\sum_{\sigma\in\Xi} \alpha_\sigma=1=\sum_{\sigma\in\Xi} \beta_\sigma$.
Then, any bijection $\Phi:{\cal M}_n(\FF)\to {\cal M}_n(\FF)$  which
maps the zeros of~$\mathfrak{p}_1$ to the zeros of $\mathfrak{p}_2$
({\small i.e., $\Phi({\Ss}_1)\subseteq {\Ss}_2$}) preserves
commutativity.
\end{theorem}

\begin{remark}
In particular if $\mathfrak{p}_1 =\mathfrak{p}_2$ the result of
Theorem~\ref{thm:main-on-vanishing-sum} also holds.
\end{remark}

\begin{corollary} \label{CorZ_2}
 Under the additional
requirement $\Phi(\Id)\not=0$, the conclusion of
Theorem~\ref{thm:main-on-vanishing-sum} is valid for~$\FF=\ZZ_2$ also.
\end{corollary}

The corollary below shows that the injectivity assumption on $\Phi$ can
be substituted by the requirement that $\Phi$ maps the zeros of
$\mathfrak{p}_1$ into the  zeros of $\mathfrak{p}_2$
 {\em strongly\/}.

\begin{corollary}\label{cor:zero-sum}%
In addition to the assumptions from
Theorem~\ref{thm:main-on-vanishing-sum}, suppose further that $k\ge 3$,
and that a ({\small possibly noninjective}) surjection~$\Phi:{\cal
M}_n(\FF)\to{\cal M}_n(\FF)$ strongly maps the zeros
of~$\mathfrak{p}_1$ to the zeros of~$\mathfrak{p}_2$.
{Then,~$\Phi(A)=0$ implies~$A=0$. Consequently,~$\Phi$ preserves
commutativity.}
\end{corollary}

The following remark provides the final forms of the transformations
satisfying conditions of Theorem~\ref{thm:main-on-vanishing-sum}.

\begin{remark}
A surjective {\em and additive} commutativity preservers are of the
form $\Phi(A)=\gamma\,TA^\varphi T^{-1}+f(A)\Id$ or
$\Phi(A)=\gamma\,T(A^\varphi)^{\tr} T^{-1}+f(A)\Id$. Here,
$\gamma\in\FF$,  $\phi:\FF\to\FF$ is a ring automorphism, and~$f:{\cal
M}_n(\FF)\to\FF$ an additive function; see Bre\v{s}ar~\cite{bresar},
Petek~\cite{petek_t}, and Beidar and Lin~\cite{beidar_lin}.

We refer to the works by \v{S}emrl~\cite{semrl:commutativity} and
Fo\v{s}ner~\cite{fosner} for a bijective, possibly non-additive
mappings, that {\em strongly preserve commutativity}. At least on the
subset of~${\cal M}_n(\CC)$, consisting of  those matrices whose Jordan
structure has only the cells of dimension at most two, they are of the
form $\Phi(A)=Tp_A(A^\varphi) T^{-1}$ or $\Phi(A)=Tp_A\bigl(
(A^\varphi)^{\tr} \bigr)T^{-1}$, where~$p_A$ is a certain polynomial
that depends on~$A$.
\end{remark}

\begin{remark}
The converse of Theorem~\ref{thm:main-on-vanishing-sum} does not hold.
Namely, there are many
polynomials~$\mathfrak{p}_1=\mathfrak{p}_2=:\mathfrak{p}$ and
commutativity preservers which {\em do not preserve the zeros
of~$\mathfrak{p}$}.
 We refer  the reader to the last section for
examples.
\end{remark}

\subsection{The proof of Theorem~\ref{thm:main-on-vanishing-sum}}
The proof will be given in a series of Lemmas.

We first recall some known results about rational forms for matrices
over an arbitrary field $\FF$.

\begin{de}
A {\em companion matrix\/} of a monic polynomial
$$f(x)=x^m+a_{m-1}x^{m-1}+\ldots+ a_3x^3+a_2x^2+a_1x+a_0,$$ of degree
$m\ge 2$, is the matrix \e \label{e0} C(f)=\left[
\begin{array}{ccccccl}
0&0&0&0&\ldots&0&-a_{0} \\
1&0&0&0&\ldots&0&-a_{1} \\
0&1&0&0&\ldots&0&-a_{2} \\
0&0&1&0&\ldots&0&-a_{3} \\
\vdots&\vdots&\vdots&\vdots&\ddots&\vdots&\vdots\\
0&0&0&0&\ldots&1&-a_{m-1} \end{array} \right]\in {\cal M}_{m}(\FF). \ee
 If $f(x)=x+a_0$ is of degree one we let $C(f):=-a_0$ be the
 $1$--by--$1$ matrix, i.e., a scalar.
\end{de}

The following lemma is straightforward and well-known:

\begin{lemma} \label{L0}
The polynomial $f$ is a characteristic polynomial of its companion
matrix $C(f)$.
\end{lemma}

\begin{theorem} \label{T1} {\rm \cite[p. 144]{Grove}, \cite[Theorem 11.20]{RML}}
Any matrix $A\in {\cal M}_n(\FF)$ is similar over $\FF$ to a matrix
${\cal C}(A)=\bigoplus_j C(p_1^{e_{1j}})\oplus \ldots \oplus
\bigoplus_j C(p_k^{e_{kj}})$, where the~$p_i$ are distinct irreducible
factors of the characteristic polynomial $\chi_A(x)=\prod\limits_{{1\le
i\le k}\atop{1\le j\le k_i}} p_i(x)^{e_{ij}}$. The matrix ${\cal C}(A)$
is determined uniquely, up to the order of diagonal blocks $C(g_i)$.
\end{theorem}

\begin{de}
The matrix ${\cal C}(A)$ described in Theorem~\ref{T1} is called a {\em
primary rational  form\/} of $A$.
\end{de}

In all statements till the end  of this section we assume that
conditions of Theorem~\ref{thm:main-on-vanishing-sum} are satisfied.

 \bigskip

\begin{lemma} \label{L1}
If~$\xi$ is nonzero scalar then the primary rational form of $\Phi(\xi
\Id)$ does not contain non-zero nilpotent blocks.
\end{lemma}
\begin{proof} Pick a similarity $P\in \mathrm{GL}_n(\FF)$ such that
$P^{-1}\Phi(\xi\Id)P$ equals the primary rational form~${\cal
C}(\Phi(\xi\Id))$ of Theorem~\ref{T1}. Now, if the claim is false, at
least one block of~${\cal C}(\Phi(\xi\Id))$ is a nonzero nilpotent. For
simplicity, assume it is  the first ({\small i.e.: the most
upper-left}) one. Therefore, it equals Eq.~(\ref{e0}), with zeros on
the last column. Then, with~$E:=PE_{11}P^{-1}$,
\e\label{e1}%
{\Phi}(\xi\Id)\,E=PE_{21}P^{-1},\quad\hbox{ and }\quad
E\,{\Phi}(\xi\Id)=0.\ee By surjectivity,  $E={\Phi}(F)$ and
${\Id}={\Phi}(J)$ for some $F,J\in {\cal M}_n(\FF)$.

Pick an integer~$t\in \{1,\ldots,k\}$ from the
Definition~\ref{def:admissible} of admissible sequence.  Note that
$t<\sigma(t)\le k$ for each $\sigma\in\Xi$. Consider now the following
matrix~$k$-tuple
$$\bigl(A_1:=J,\dots,A_{t-1}:=J,
{\mathbf A}_{t}:={\mathbf {(\xi\Id)}},A_{t+1}:=F,\dots,A_{k}:=F
\bigr),$$
 which lies in~${\Ss}_1$, since each
$\sigma$ fixes the indices $\{1,\ldots,t-1\}$, and
$\sum\alpha_\sigma=1$. By the
assumptions,~$\Phi({\Ss}_1)\subseteq{\Ss}_2$, and we have
$$ \Phi(A_1)\cdots \Phi(A_k)=\sum_{\sigma\in\Xi} \bigl(\beta_\sigma
\Phi(A_{\sigma(1)})\cdots
\Phi(A_{\sigma(k)})\bigr),$$
 i.e., \e \label{e2}
\Id^{t-1}\Phi(\xi\Id)E^{k-t}=\sum_{\sigma\in\Xi}
(\beta_\sigma\Id^{t-1}E^{g_\sigma}\,\Phi(\xi\Id)E^{k-t-g_\sigma}), \ee
 where
$g_\sigma:=\sigma(t)-t>0$. However, the matrix $E$ is idempotent, so
$E^{k-t}=E=E^{g_\sigma}$, and it follows from~(\ref{e1}) that the left
hand side of the equality~(\ref{e2}) is equal to
${\Phi}(\xi\Id)\,E=PE_{21}P^{-1}$, while the right hand side is equal
to 0, a contradiction. \end{proof}

\begin{lemma} \label{L2}
Suppose~$\FF\not=\ZZ_2$. Then, there exists a nonzero scalar~$\xi$ such
that $\Phi(\xi\Id)\in \mathrm{GL}_n(\FF)$.
\end{lemma}
\begin{proof}Since~$\Phi$ is injective and the cardinality
$|\FF\backslash\{0\}|\ge 2$, there
exists at least one nonzero scalar~$\xi$ such
that~$\Phi(\xi\Id)\not=0$. As in the proof of Lemma~\ref{L1}, let~$t$
be fixed by Definition~\ref{def:admissible} let
$g_\sigma:=\sigma(t)-t>0$, and let $J:=\Phi^{-1}(\Id)\in {\cal
M}_n(\FF)$. Here, we consider the following matrix $k$-tuple:
$$\bigl(A_1:=J,\dots,A_{t-1}:=J,
{\mathbf A}_{t}:={\mathbf X},A_{t+1}:={(\xi \Id)},\dots,A_{k}:={(\xi
\Id)} \bigr).$$
 Again, this $k$-tuple is in ${\Ss}_1$ for an arbitrary
matrix~$X$. By the assumptions, $\Phi({\Ss}_1)\subseteq{\Ss}_2$, and we
have
\begin{equation}\label{eq:Phi(Id)-is-nonsingular}%
 \Id^{t-1}\Phi(X)\Phi(\xi \Id)^{k-t}=
 \sum_{\sigma\in\Xi}\bigl(\beta_\sigma \Id^{t-1}\Phi(\xi \Id)^{g_\sigma}
  \Phi(X)\Phi(\xi \Id)^{k-t-g_\sigma}\bigr),
\end{equation}
Let us assume that $\Phi(\xi \Id)$ is singular.
Recall~$\Phi(\xi\Id)\not=0$, so by Lemma~\ref{L1}, the primary rational
form, ${\cal C}(\Phi(\xi \Id))$, contains at least one zero block and
at least one non-zero block.  For simplicity, assume the first one is
zero, i.e., ${\cal C}(\Phi(\xi \Id))=\mathbf{0}\oplus \mathbf{C}$,
where~$\mathbf{C}\not=0$ is a sum of all, but the first, blocks.

By Lemma~\ref{L1}, ${\cal C}(\Phi(\xi \Id))^{k-t}=\mathbf{0}\oplus
\mathbf{C}^{k-t}\not=0$. Consequently, the matrix ${\cal C}(\Phi(\xi
\Id))^{k-t}$ has a nonzero row, i.e., there exists $p$, $1\le p\le n$
such that {$E_{1p}\,{\cal C}(\Phi(\xi \Id))^{k-t}\not=0$}. However,
note that ${\cal C}(\Phi(\xi \Id))E_{1p}=0$, so also ${\cal C}(\Phi(\xi
\Id))^{g}E_{1p}=0$ for all $g\in\NN\backslash\{0\}$, in particular for
each $g:=g_\sigma$. Now, consider $P\in \mathrm{GL}_n(\FF)$ such that $
\Phi(\xi \Id)=P{\cal C}(\Phi(\xi \Id))P^{-1}$, and choose a matrix~$X$
with $\Phi(X)=PE_{1p}P^{-1}$. For such~$X$, the left hand side
of~(\ref{eq:Phi(Id)-is-nonsingular}) is nonzero while the right hand
side is zero, a contradiction. \end{proof}

\begin{lemma} \label{L3}
$\Phi$ preserves commutativity.
\end{lemma}
\begin{proof} Without loss of generality we assume that $\xi=1$ in
Lemma~\ref{L2}, i.e., that~$\Phi(\Id)\in \mathrm{GL}_n(\FF)$ ---
otherwise, the bijection~$\Phi(\xi\cdot \anyargument)$ would be
considered instead of~$\Phi$.
 By the definition of admissible sequence, there exists
 an
integer $w$ such that $u\equiv\sigma(w+1) < \sigma(w)\equiv
v\;\forall\,\sigma\in\Xi$, we fix the smallest such index $w$.
 Consider the following matrix $k$-tuple: \e
\label{e3} (A_1:=\Id,\ldots,A_{w-1}:=\Id, {\mathbf A}_w:={\mathbf X},
{\mathbf A}_{w+1}:={\mathbf Y}, A_{w+2}:=\Id,\ldots,A_k:=\Id). \ee
 If $XY=YX$
then the $k$-tuple~(\ref{e3}) is in ${\Ss}_1$. Thus,
\begin{equation}\label{e4}%
\begin{split}
\Phi(\Id)^{w-1}\,\Phi(X)\Phi(Y)\,&\Phi(\Id)^{k-w-1}= \\
&=\sum_{\sigma\in\Xi}\beta_\sigma \Phi(\Id)^{u-1}\,\Phi(Y)\,
\Phi(\Id)^{s}\,\Phi(X)\,\Phi(\Id)^{k-v},
\end{split}\end{equation}
where ${s}=v-u-1\ge 0$. Note that by the definition of admissible
sequence $u$ and $v$ are independent of~$\sigma$ and thus the right
hand side is equal to
$\Phi(\Id)^{u-1}\,\Phi(Y)\,\Phi(\Id)^{s}\,\Phi(X)\,\Phi(\Id)^{k-v}$.
Since~$\Phi(\Id)$ is invertible, Eq.~(\ref{e4}) simplifies into:
\begin{equation}\label{e4.5}%
\Phi(X)\Phi(Y)= \Phi(\Id)^{u-w}\,\Phi(Y)\,\Phi (\Id)^s\,
\Phi(X)\,\Phi(\Id)^{w-v+1}
\end{equation}
whenever $XY=YX$.
 We
first claim that~$\Phi(\Id)^{s}$ is a scalar matrix. Assume on the
contrary, and consider two cases:\medskip

{\bf Case 1}: The primary rational form~${\cal
C}(\Phi(\Id)^{s})=P^{-1}\Phi(\Id)^{s}P$ contains a block of
dimension~$\ge 2$. Again, for the sake of simplicity,  we assume this
is the first block. Therefore, the~$(1,1)$-entry of~${\cal
C}(\Phi(\Id)^{s})$ is zero. We would then let~$X=Y$ be such
that~$\Phi(X)=PE_{11}P^{-1}=\Phi(Y)$. This contradicts~(\ref{e4.5}),
since the left hand side would be $PE_{11}P^{-1}$, while the right
would be zero.\medskip

{\bf Case 2}: All blocks of ${\cal
C}(\Phi(\Id)^{s})=P^{-1}\Phi(\Id)^{s}P$ are one-dimensional, i.e.,
${\cal C}(\Phi(\Id)^{s})=\diag(d_1,\dots,d_n)$ is diagonal. Since
${\cal C}(\Phi(\Id)^{s})$ is non-scalar, at least two diagonal entries
differ. For the sake of simplicity, assume~$d_1\not=d_2$. Note that by
Lemma~\ref{L2} the matrix~$\Phi(\Id)^{s}$ is invertible,
so~$d_1\not=0$. Then, the
similarity by the matrix $S:=\left[\begin{smallmatrix} \frac{d_2}{d_1} & 1\\
1 & 1
\end{smallmatrix}\right]\oplus \Id_{n-2}$ transforms this matrix into

$$S \,{\cal C}(\Phi(\Id)^{s})\,S^{-1}=\begin{bmatrix}
0&d_2&\\
-d_1&d_1+d_2 &
\end{bmatrix}\oplus\diag(d_3,\dots,d_n)$$
with the zero $(1,1)$-entry. We can then reach a contradiction as in
Case~1.

Consequently, since~$\Phi(\Id)$ is invertible we
have~$\Phi(\Id)^s=\lambda \Id\not=0$, and Eq.~(\ref{e4.5}) further
simplifies  into:
\begin{equation}\label{eq:pom1}%
\Phi(X)\Phi(Y)=\lambda
\Phi(\Id)^{u-w}\,\Phi(Y)\Phi(X)\,\Phi(\Id)^{w-v+1}
\end{equation}
whenever $XY=YX$. Pick any distinct indices~$i,j\in\{1,\dots,n\}$.
Since~$n\ge 3$, there exists another one,
$k\in\{1,\dots,n\}\backslash\{i,j\}$. Now, by the surjectivity, we may
choose~$X=Y$ such that $\Phi(X)=\Phi(Y)=E_{ik}+E_{kj}$, which gives
$\Phi(X)\Phi(Y)=(E_{ik}+E_{kj})^2=E_{ij}=\Phi(Y)\Phi(X)$. We can
likewise find~$X=Y$ such that
$\Phi(X)\Phi(Y)=E_{ii}^2=E_{ii}=\Phi(Y)\Phi(X)$. Putting this into
Eq.~(\ref{eq:pom1}) we deduce that $E_{ij}=\lambda
\Phi(\Id)^{u-v}\,E_{ij}\,\Phi(\Id)^{w-v+1}$ for any indices~$(i,j)$.
Consequently,
$$A=\lambda
\Phi(\Id)^{u-v}\,A\,\Phi(\Id)^{w-v+1};\qquad \forall\,A\in{\cal
M}_n(\FF).$$
 With~$A:=\Id$ we have $\Phi(\Id)^{u-v}\cdot
\Phi(\Id)^{w-v+1}=1/\lambda \,\Id$. Therefore, $AS=(\lambda S)A$ for
each~$A$, where $S:=\Phi(\Id)^{u-v}$. A standard procedure
with~$A:=E_{ij}$ gives that~$S$ is scalar, and~$\lambda=1$. {Hence,
$\Phi(\Id)^{w-v+1}=S^{-1}$ is also a scalar matrix}. Since $\lambda=1$,
this further reduces Eq.~(\ref{eq:pom1}) into the desired
$\Phi(X)\Phi(Y)=\Phi(Y)\Phi(X)$ whenever~$XY=YX$.
 \end{proof}\bigskip

\subsection{Proof of Corollaries}

\begin{remark}\label{RKer}%
As it was seen in the proof of Theorem~\ref{thm:main-on-vanishing-sum},
we used the requirement of injectivity just once: in the proof of
Lemma~\ref{L2}. {Moreover, even there it suffices to have its curtailed
form, i.e.,  that} there exists $\xi\in \FF$, $\xi\ne 0$, such that
$\Phi(\xi \Id)\ne 0$. Thus the result of
Theorem~\ref{thm:main-on-vanishing-sum} holds under these, even more
general, conditions.
\end{remark}
\begin{proof}[Proof of Corollary~\ref{CorZ_2}]
{By the assumptions, $\Phi(\Id)\ne 0$, so} the proof of Lemma~\ref{L2}
works by Remark~\ref{RKer}. It can be directly checked that the rest of
the proof of Theorem~\ref{thm:main-on-vanishing-sum} does not use the
assumptions on the ground field.
\end{proof}

\begin{proof}[Proof of Corollary~\ref{cor:zero-sum}]
{Suppose $\Phi(A)=0$ for some matrix~$A$, and} consider the following
matrix $k$-tuple
$$\bigl(X_1=\Id,\dots,X_{w-1}=\Id,{\bf X}_{w}={\bf A},{\bf
X}_{w+1}={\bf X},X_{w+2}=\Id,\dots,X_k=\Id\bigr).$$
 This $k$-tuple is mapped
by~$\Phi$  into a $k$-tuple with~$w$-th member equal to
$\Phi(X_w)=\Phi(A)=0$. Therefore,
$(\Phi(X_1),\dots,\Phi(X_k))\in{\Ss}_2$, so that also
$(X_1,\dots,X_k)=(\Id,\dots,\Id,X_w=A,X_{w+1}=X,\Id,\dots,\Id)\in{\Ss}_1$,
for every choice of~$X$. Since~$\sigma(w+1)<\sigma(w)$ for every
admissible permutation, this further yields
$$\Id^{w-1}AX\Id^{k-w-1}=\sum_\sigma \alpha_\sigma\Id^i X\Id^j
A\Id^l$$ for some integers~$i,j,l$. This immediately simplifies into
$AX=\sum\limits_\sigma \alpha_\sigma XA=XA$ for every~$X$. Hence,~$A$
has to be scalar.\medskip

It remains to show~$A=0$. Consider now the following $k$-tuple:
$X_{w}:=Y,X_{w+1}:=Z$, and $X_i:=A$ for the rest of indices ({\small
since~$k\ge 3$, at least one member equals~$A$}). As before we deduce
that this $k$-tuples is in~${\Ss}_1$ for every choice of~$Y,Z$. Hence,
since~$A$ is scalar, and $\sigma(w+1)<\sigma(w)$,
$$A^{k-2} YZ=A^{k-2}\sum_{\sigma\in\Xi} \alpha_\sigma ZY=A^{k-2}ZY$$
for any choice of~$Y,Z$. This is possible only when~$A^{k-2}=0$,
i.e.,~$A=0$.\medskip

Therefore, nonzero scalar matrices are not annihilated by~$\Phi$, and
we can follow the proof of Theorem~\ref{thm:main-on-vanishing-sum} to
see that~$\Phi$ preserves commutativity.
\end{proof}

\section{\label{S5} Concluding remarks and examples}
This section mainly contains counterexamples  to show that our results
cannot be further improved without imposing additional hypothesis.

{\em The following example shows that the inverse implication does not
hold in Theorem~\ref{thm:main-on-vanishing-sum}, namely there are
commutativity preserving mappings that do not preserve zeros of a fixed
polynomial.\/}

\begin{example}
Let $\Phi(A)=A+a_{12}\,{\rm Id}$ for all $A=[a_{ij}]\in {\cal
M}_n(\FF)$. We consider the polynomial $\mathfrak{p}(x,y,z):=xyz-yxz$.
Then the triple $(x:=E_{11},\, y=z:=E_{12})$ is a zero of this
polynomial, but its image $(E_{11},\Id+E_{12},\Id+E_{12})$ is not a
zero of this polynomial.
\end{example}

{\em There exist ({\em even linear}!) transformation $\Phi$ which
strongly preserve commutativity, but have a nonzero kernel.\/}
\begin{example}
 Let us consider $\Phi:{\cal M}_n(\FF) \to {\cal M}_n(\FF)$ defined by
$\Phi(A):=A-\mathrm{Tr}(A)/n\, \Id$. Then~$\Phi$ is linear, $\Phi$
strongly preserves commutativity, but $\Phi(\Id)=0$.
\end{example}

{\em Corollary~\ref{cor:TAllMatr} is no longer true when~$k=2$, i.e.,
the transformation $\Phi$ may not be controllable on some large subsets
of $\DPhi_1$.\/}

\begin{example}\label{exa:counter-axa-for-nonvanishing}
 Namely,
consider $\mathfrak{p}(x,y):=xy+yx$. Pick any  $A\in\Theta:=\{A\in
 {\cal  M}_n(\FF);\;\;0\not\in\Sp(A)+\Sp(A)\}$. Then, by
 Sylvester--Rosenblum Theorem
({\small cf.\ the proof of Lemma~\ref{lem:Sylvester--Rosenblum}}), the
elementary operator $X\mapsto XA+AX$ is invertible, so
$\mathfrak{p}(X,A)=0$ if and only if~$X=0$. Hence, the restriction
of~$\Phi$ to the subset~$\Theta$
 has no structure at all, i.e.,~$\Phi|_{\Theta}$ may arbitrarily permute
the elements of~$\Theta$, yet it still strongly preserves the zeros
of~$\mathfrak{p}(x,y)$. We remark that~$\Theta$ is a rather large
subset: when~$\FF=\CC$ it is nonempty and open in~${\cal M}_n(\CC)$, by
continuity of the eigenvalues.
\end{example}

{\em The characterization of
Lemma~\ref{lem:charact_of_minimal_idempotents} is no longer valid
if~$\FF$ is not algebraically closed.\/}

\begin{example}
 As a
counterexample,  choose   a polynomial~$\mathfrak{p}(x,y):=xy+yx$
 and consider $A:=\left[\begin{smallmatrix}
 0 & -3 \\
 1 & 0
\end{smallmatrix}\right]\oplus 1\in{\cal  M}_3(\IR)$.
Then,~$\Omega_{\bullet A}=\Omega_{A\bullet A}=\{\left[\begin{smallmatrix}
 -d & 3 c  \\
 c & d
\end{smallmatrix}\right]\oplus 0;\;\;c,d\in\IR\}$, and it contains no
nontrivial square-zero matrices, nor it equals~$\IR\,P$ for
idempotent~$P$.
\end{example}

{\em The following example shows that there are field automorphism
which do not preserve zeros of matrix polynomials.\/}
\begin{example}
\label{ExPhi} Let us consider the polynomial $\mathfrak{p}=xy-i\,yx$
and  automorphism $\varphi:{\mathbb C}\to {\mathbb C}$ which sends any
$x\in {\mathbb C}$ to its complex conjugated element $\overline x$.
Then consider $$ A:=\begin{bmatrix}
 1 & 1 \\ -1 & 1 \end{bmatrix}\oplus \Id_{n-2}\in {\cal M}_n(\mathbb C),\
B:=\begin{bmatrix} 1 & i \\ i & -1\end{bmatrix}\oplus \Id_{n-2}\in
{\cal M}_n(\mathbb C).$$ The direct computations show that the matrices
$A,B$ are zeros of $\mathfrak{p}$ but their conjugated matrices
$A^{\varphi}, B^{\varphi}$ are not.
\end{example}

{\em Finally, the transposition transformation may not preserve zeros
of $\mathfrak{p}$ for some $\mathfrak{p}$.\/}
\begin{example}
\label{Ex_tr} Let us consider the polynomial $\mathfrak{p}=xy$ and the
matrices $A=E_{11}$, $B=E_{21}$. Then $\mathfrak{p}(A,B)=0$, however
$\mathfrak{p}(A^\tr,B^\tr)=E_{11}E_{12}=E_{12}\ne 0$.
\end{example}

{\bf Acknowledgments}. The authors are indebted to Professor Chi-Kwong
Li for inspiring conversations regarding the topic of Chapter~\ref{S4}.

\end{document}